\documentclass{ieeeaccess}
\usepackage{cite}
\usepackage{algorithmic}
\usepackage{graphicx}
\usepackage{textcomp}
\def\BibTeX{{\rm B\kern-.05em{\sc i\kern-.025em b}\kern-.08em
    T\kern-.1667em\lower.7ex\hbox{E}\kern-.125emX}}
\usepackage{booktabs}
\usepackage{dcolumn}
\usepackage{bm}
\usepackage{xcolor}
\usepackage{hyperref}
\usepackage{tabularray}
\usepackage{multirow}
\usepackage{amssymb}
\usepackage{amsfonts}
\usepackage{amsmath}
\usepackage{mathtools}
\usepackage{enumitem}
\usepackage{tabularx}
\usepackage[font=scriptsize]{caption}
\usepackage{subfigure}
\usepackage{adjustbox}

\renewcommand{\thetable}{\Roman{table}}

\usepackage[classfont=sanserif,
langfont=sanserif,
funcfont=sanserif ]{complexity} 
\newclass{\KP}{KnapsackProb}
\newclass{\CO}{ColOpt}
\newclass{\XORSAT}{XORSAT}
\newclass{\NISQ}{NISQ}
\newclass{\TSP}{TravSalProb}
\newclass{\BinPP}{BinPackProb}
\newclass{\MIS}{MaxIndSet}
\newclass{\MAXCUT}{MAXCUT}

\usepackage{array}
\usepackage{makecell}

\DeclareMathOperator{\bin}{bin}
\DeclareMathOperator{\dec}{dec}

\begin{document}
\history{Date of publication xxxx 00, 0000, date of current version xxxx 00, 0000.}
\doi{10.1109/TQE.2020.DOI}

\title{
Approaching Collateral Optimization for NISQ and Quantum-Inspired Computing (May 2023)
}
\author{\uppercase{Megan C. Giron\authorrefmark{1},
 Georgios Korpas\authorrefmark{1,}\authorrefmark{2},Waqas Parvaiz\authorrefmark{1}, Prashant Malik\authorrefmark{3} and Johannes Aspman}.\authorrefmark{2}}

\address[1]{ HSBC Lab, Innovation \& Ventures,  8 Canada Square, London, E14 5HQ, U.K.}
\address[2]{Department of Computer Science, Czech Tech. University in Prague, Karlovo nam. 13, Prague 2, Czech Republic}
\address[3]{Markets and Security Services,  8 Canada Square, London, E14 5HQ, U.K.}


\markboth
{Giron \headeretal: Approaching Collateral Optimization for NISQ and Quantum-Inspired Computing}
{Giron \headeretal: Approaching Collateral Optimization for NISQ and Quantum-Inspired Computing}

\corresp{Corresponding authors: Megan C. Giron; Georgios Korpas; Waqas Parvaiz (email:megan.giron@hsbc.com; georgios.korpas@hsbc.com; waqasparvaiz@live.co.uk).}

\begin{abstract}
Collateral optimization refers to the systematic allocation of financial
assets to satisfy obligations or secure transactions, while simultaneously
minimizing costs and optimizing the usage of available resources. {This involves assessing number of characteristics, such as cost of funding and quality of the underlying assets to ascertain the optimal collateral quantity to be posted to cover exposure arising from a given transaction or a set of transactions. One of the common objectives is to minimise the cost of collateral required to mitigate the risk associated with a particular transaction or a portfolio of transactions while ensuring sufficient protection for the involved parties}. Often, this results in a large-scale combinatorial optimization problem. In this study, we initially present a Mixed Integer Linear Programming (MILP) formulation for the collateral optimization problem, followed by a Quadratic Unconstrained Binary optimization (QUBO) formulation in order to pave the way towards approaching the problem in a hybrid-quantum and NISQ-ready way. We conduct local computational small-scale tests using various Software Development Kits (SDKs) and discuss the behavior of our formulations as well as the
potential for performance enhancements. We find that while the QUBO based approaches fail to find the global optima in the small scale experiments, they are reasonably close suggesting their potential for large instances. We further survey the recent literature that proposes alternative ways to attack combinatorial optimization problems suitable for collateral optimization.
\end{abstract}


\titlepgskip=-15pt

\maketitle

\section{Introduction}
\label{sec:level1}
\PARstart{I}{n} In the context of a financial transaction, wherein one party lends assets to another, the lender assumes a credit risk arising from the possibility that the counterparty may default on their obligations. This risk also arises in derivatives transactions where the party ``in-the-money'' is exposed to the party ``out-of-the-money''. To mitigate this risk, the borrower is required to provide low-risk securities (such as cash, bonds, or equities) to the lender for the duration of the transaction. This practice, known as \emph{collateralization} \cite{Simmons2019}, serves as a form of security against loan defaults as the lender can seize these assets to offset any losses resulting from default. The value of the collateral received is expected to be commensurate with the outstanding exposure, in order to effectively counterbalance the associated risk.

Often, a bilateral {contract (or schedule)} is formed to agree on the terms under which securities can be considered collateral, the process of evaluating the value of these assets, and other regulations. The relevant party may then accordingly select the assets they post to the counter-party. For large financial institutions, such as banks, this can involve a pool of numerous assets to choose from which need to be distributed amongst a portfolio of various counterparties (other banks, hedge funds, central banks, etc.). Each asset has an associated opportunity cost, which is a measure of how valuable the asset would be if it were used for another purpose. As well as a cost related to the risk of posting to a particular counterparty amongst other administrative costs. The bank must therefore carefully consider their choice of transactions to reduce these costs. However, the magnitude of the possible combinations of allocations for large institutions makes this a time-consuming process. In addition, poor collateral management can have significant consequences. The 2008 financial crisis was partly due to collateralization of high-risk securities and led to the bankruptcy of some of the largest financial institutions, among other consequences \cite{N_tzenadel_2020}. 

This crisis has led to the reformation of many financial processes through the implementation of regulations such as Basel III \cite{BaselIII}, Dodd-Frank Act \cite{Dodd}, EMIR \cite{EMIR} as well as motivating academic research with regards to how collateral can be better managed \cite{brigo2011counterparty,deloitte2014}. These studies relating to collateral management are generally centered on financial theories such as risk aversion and its global financial impact. A crucial aspect of collateral management is the development of an automated process that selects the optimal allocations. Despite the importance of collateral optimization (\CO), the literature surrounding this topic is sparse due to the competitive advantage these strategies offer to financial institutions. 

Naturally, linear programming algorithms can provide a framework for tackling such problems \cite{Bertsimas1997-pw}. Specifically, \CO~is suitable to be implemented using mixed-integer linear program (MILP) solvers such as the ones available with IBM CPLEX \cite{cplex2009v12}, Gurobi \cite{gurobi}, or Mosek \cite{mosek}. The success of a given \CO~instance, or the quality of its solution, is dependent on a clear mapping between the business problem and the mathematical formulation as well as the choice of the precise implementing algorithm. The benefit of using numerical optimization is that we attain the allocation selections in a single process, which is in contrast to other proposed models such as ``ranking-based", ``economic-cost" and ``waterfall" models which are sequential in nature, rather than automated \cite{bylund2017collateral}. 
However, there are several limitations of MILP solvers when applied to problems such as \CO, which potentially involve complex nonlinearities and large-scale datasets. For example, MILP solvers have exponential worst-case complexity and can take a significant amount of time to solve large-scale complex optimization problems. While having convergence certificates is very desirable, the exponential complexity is often a problem for many \CO~instances that involve a large number of decision variables and constraints and it is not uncommon to either have long solution times or even infeasibility. That said, MILP solvers are the standard in both industrial and academic applications but the community is very keen on exploring alternative approaches. 

Different avenues to (approximately) solve such computationally challenging problems could be provided through 
alternative computational models. For example, in Ref. \cite{Mniszewski2019}, IBM's True Spike neuromorphic computer \cite{Schuman2022} was utilized to find approximate solutions for the graph partitioning problem. More interest is directed towards quantum computing \cite{ladd2010quantum}. These computers rely on quantum-mechanical effects for storing and processing information. However, because of the fragile environment required for these effects to occur, the realization of fault-tolerant quantum computers is still a difficult task.  

Despite this, it is believed that Noisy Intermediate-Scale Quantum (NISQ) era devices can provide an advantage in the finance industry since these business use cases can be well formulated for near-term quantum devices \cite{Egger2020, herman2022survey}. The field of ``quantum finance'' can be divided into three sections: stochastic modeling (for example, quantum alternatives to Monte Carlo simulations) \cite{egger2020quantum,bouland2020prospects,herman2022survey, rebentrost2022quantum,intallura2023survey}, machine learning \cite{egger2020quantum,bouland2020prospects,pistoia2021quantum,herman2022survey,jacquier2022quantum,emmanoulopoulos2022quantum,kirk2023emergent} and optimization \cite{rosenberg2015solving,egger2020quantum,rebentrost2018quantum,leclerc2022financial,lim2022quantum,jacquier2022quantum,brandhofer2023benchmarking,vesely2023finding,sakuler2023real}, all of which have had a recent gain in interest. 

{
Very often, in this context, the prototypical optimization use-case is that of (Markowitz) portfolio optimization. While this use-case shares a few similarities with \CO, there are a few fundamental differences as well. The constraints of \CO, seem to be more involved even within the simplifications we present in Sec. \ref{sec:level2}. An important difference is also the fact that our objective function is inherently linear, at first glance. However, several ``tricks'' can be performed in order to end up with a well-behaved formulation suitable for a variety of Ising NISQ or hybrid solvers. 
}

Quantumly, there are many approaches to follow in order to approximate better solutions for a variety of \NP-Hard problem instances encountered in finance. The main three approaches (listed below)  have a common theme: conversion of the original mathematical formulation of the problem from a linear program (LP) formulation to a quadratic unconstrained binary optimization problem (QUBO). The reason lies in the inherent ability of quantum or hybrid approaches by modeling the Ising model type of systems (see App. \ref{ap:qubo_to_ising}). There, an optimization problem is mapped to the classical
Hamiltonian of the Ising model, where its ground state encodes the optimum. As a matter of fact, many \NP-Hard problems, including Karp's list of 21 \NP-Hard problems, are known to admit at least one formulation of the Ising model \cite{lucas2014ising}. The QUBO or Ising approach can be used and problems can be mainly tackled as follows:
\begin{itemize}
    \item Using variational quantum algorithms (VQAs) \cite{cerezo2021variational} such as the Quantum Approximate Optimization Algorithm (QAOA) \cite{farhi2014quantum} on gate-based quantum computers (e.g., IBM's superconducting quantum computer). A variety of tests have been performed in this context with significant (qubit) resource improvements recently \cite{chatterjee2023solving}. 
    \item Using quantum annealing (QA) \cite{Apolloni1989,Finnila1994,Farhi2001,RevModPhys.90.015002,hauke2020perspectives} on adiabatic quantum computers (quantum annealers) such as the D-Wave hardware \cite{willsch2022benchmarking} (see \cite{phillipson2021portfolio} for an application to portfolio optimization).
    \item Using quantum-inspired methods which can be understood as using the QUBO formulation of the problem of interest with any approach that ranges from simulated annealing \cite{Van_Laarhoven1987-kd} on high-performance clusters or digital annealing \cite{aramon2019physics} such as Fujitsu's FPGA-based ``quantum-inspired'' classical hardware DAU \cite{fujitsu}.
\end{itemize}

These approaches are very promising and it is widely expected that near-term quantum computers, as well as the dedicated quantum inspired hardware, may have a good chance to provide computational or business advantage\footnote{By the term ``business advantage'' we mean enough computational resource savings or improvement in performance enough to justify the adoption of the underlying technology.} \cite{PhysRevApplied.19.024027,lubinski2023optimization} in the near term. Note that the formulation we will present below in Sec. \ref{sec:level2}, are hardware agnostic and in that sense any of the above approaches would be suitable in principle.

Additionally, there exists evidence that the class of problems that NISQ computers can solve is not a subset \footnote{Ref. \cite{chen2022complexity} defines \NISQ~as the class that contains all problems that can be solved
by a polynomial-time probabilistic classical algorithm with access to a noisy quantum device and where \BPP$\subseteq$\NISQ$\subseteq$\BQP.} of \BPP \cite{chen2022complexity}.
However, we realize that their heuristic nature can conceal their applicability, especially when considering problem instances with sizes suitable for high-quality MILP solvers. Such a result was reported in Ref. \cite{pusey2020adiabatic} wherein the authors discuss how D-Wave's  2000Q machine (as well as classical simulated annealing) failed to even come close to the branch and bound approaches \cite{Papadimitriou,huang2021branch} for solving certain instances of the Knapsack Problem (\KP)~(see Sec. \ref{sec:knapsack}). One can potentially try to use the VQA approach instead \cite[for example]{cerezo2021variational, bouland2020prospects, herrman2021globally,Amaro_2022}. In this context, Nannicini \cite{PhysRevE.99.013304} reported that there seems to be a lack of transparent computational advantage in introducing entanglement when VQAs are used. This lack of computational advantage is expected to some extent due to the well-known local minima problem that VQAs exhibit \cite{bittel2021training} (see \cite{rivera2021avoiding} for the proposed way to avoid this problem). In addition, bias in the noise of circuits that implement VQAs can unfavourably affect the convergence ratios \cite{kungurtsev2022iteration}. Other limitations are discussed in \cite{PRXQuantum.4.010309}.

However, performance advantages have been showcased in a variety of applications within the context of digital and quantum annealing solutions, tensor networks, and analog and digital (gate-based) quantum computing. For example, Ref. \cite{ebadi2022quantum} studied the maximal independent set (\MIS) problem and found a superlinear quantum speedup, as opposed to classical solutions, when considering very hard graphs. From the point of view of computational complexity, the \MIS~problem is not particularly different from other \NP-Hard problems such as the \KP~(which actually is weakly \NP-Complete) or other \NP-Hard problems that admit a MILP formulation, and as such the results of \cite{ebadi2022quantum} are encouraging for other problems as well. An interesting benchmark test \cite{tasseff2022emerging} using the D-Wave machine showed very promising results. It is worth mentioning that D-Wave recently announced \cite{King2023} the largest quantum simulation done, in different context (spin glasses) to what is of interest here. 

In this paper, we study the \CO~problem in detail and provide a MILP formulation that we use as a testbed for: \emph{(i)} formulating a QUBO version of \CO~which makes it suitable for feeding into quantum and quantum-inspired solvers and \emph{(ii)} performing small scale simulations of such solvers and comparing to MILP. Specifically, in our MILP formulation we choose our objective to be minimization of the cost of posting collateral (different approaches to the objective of \CO~have been proposed, for example, see \cite{bylund2017collateral}). We survey and try numerous QUBO encodings and find that, modulo emulator limitations, in such small instances the quantum-inspired formulations perform well enough to be promising for implementation on real quantum or quantum-inspired hardware for very large instances. The structure of the paper is as follows. In Sec. \ref{sec:knapsack} we provide an overview of some of the different QUBO formulations for the \KP~problem. This section serves both as an introduction to the concepts used throughout the paper and to inform the \CO~problem that follows. In Sec. \ref{sec:level2} we provide the MILP formulation of the \CO~problem as well as a few QUBO proposals. In Sec. \ref{sec:illustrations} we provide a few numerical results using the formulation of the previous section. Finally, we summarize and conclude in Sec. \ref{sec:summary}.


We want to clarify that while our paper investigates various QUBO formulations for the \KP~problem, our ultimate goal is to apply this information to the \CO~problem. Specifically, we plan to use the best-performing QUBO formulations from our \KP~study to formulate and solve the collateral optimization problem using QUBO. However, we would like to note that our paper does not aim to provide an empirical comparison between quantum and classical approaches for solving MILPs, given the limited computational resources available to us (see \cite{Jnger2021} for work on the comparison of classical and quantum (adiabatic) optimization, wherein the authors discovered ``surprising'' results favoring the D-Wave machine). Additionally, we focus on small problem instances only, and we acknowledge that all the results presented in our study are heuristic in nature. Nevertheless, based on the literature results on the potential of QUBO formulations (quantumly and not only), we believe that the formulations we propose may have value in tackling larger instances of the collateral optimization problem and therefore may warrant further investigation. 

In summary, the main objective of our paper is to present a case study on the formulation and approach of the \CO~problem using quantum computing techniques, with the overarching aim of advancing the ongoing effort towards achieving ``quantum advantage'' in practical applications.


\section{Interlude with the Knapsack Problem}\label{sec:knapsack}
To inform and ensure that our formulations and the subsequent computations, we perform a simple test using a small \KP~instance. In essence, \KP~
involves determining the optimal approach to filling a knapsack of capacity $W$ with the highest possible value from a set of $n$ items that have specific sizes and corresponding values (see Table \ref{tab:KPinstance}). This problem is of interest due to its simplicity to formulate, and its simple constraints. For us, it is further interesting since we view the \CO~problem as a (somewhat complicated) generalization.

Given the large number of ``hard'' constraints of the \CO~problem, see Sec. \ref{sec:level2}, we aim to compare formulation for small instances of the Knapsack Problem with the hope to inform our approach for collateral optimization. The ``standard route'' to encode constraints to a QUBO model is by using (balanced) slack variables for penalization \cite{glover2018tutorial} (see App. \ref{ap:qubo_to_ising}). A different approach is using ``unbalanced'' penalization \cite{montanez2022unbalanced} which turns out to be particularly useful for QAOA solutions as it reduces the resources required by a gate-based machine. 

The MILP formulation of the \KP~is a well-known and straightforward approach. In the problem instance we consider, we are given a set of weights $w\in \mathbb{Z}_{\geq 0}^n$ and their corresponding values $v \in \mathbb{Z}_{\geq 0}^n$, and the objective is to maximize the total value of the items that can be packed into a knapsack subject to a given weight limit. The problem can be mathematically defined as follows:
\begin{equation}\label{eq:kp}
    \begin{aligned}
        \max_{} & \quad \sum_{i=1}^n v_i x_i, \\
        {\text{s.t.}} &\quad \sum_{i=1}^n w_i x_i \leq W,
    \end{aligned}
\end{equation}
where $W$ is the maximum weight limit (threshold) of the knapsack and $x_i$ is the binary variable representing whether the $i$-th item is to be placed in the knapsack. The best running-time algorithm for solving the \KP~is based on dynamic programming with pseudo-polynomial complexity $\mathcal{O}(d_nW)$ \cite{axiotis2018capacitated}, where $d_n$ is the number of distinct weights available while near-linear running times algorithms (in $d_n,W$) were documented  in \cite{bateni2018fast}.


The specific problem instance we consider in our study, as outlined in \cite{Silvano}, comprises ten items and possesses a known optimal solution. We leverage this knowledge to heuristically evaluate the effectiveness of our approach and guide our efforts towards tackling the larger \CO~problem discussed in Section \ref{sec:level2}. The relevant (toy) data pertaining to the items in this problem instance can be found in Table \ref{tab:KPinstance}.

\begin{table}[]
    \centering
    \begin{tabular}{l|cccccccccc} \toprule[2pt]
       Object label & A&B&C&D&E&F&G&H&I&J \\ \hline
       
       Weight& 23&31&29&44&53&38&63&85&89&82 \\ 
        Value& 92&57&49&68&60&43&67&84&87&72 \\ \hline
    \end{tabular}
    \caption{Input data used of the specific \KP~instance considered in this paper. The total capacity of the knapsack is 165.}
    \label{tab:KPinstance}
\end{table}

Although the \KP~is known to be (weakly) \NP-complete, simple instances such as the one we consider in this study can be efficiently solved by a range of classical solvers. For our experimental analysis, we used the HiGHS \cite{Huangfu2017} and GLPK \cite{GLPK} solvers, both of which yielded solutions that were in agreement with the known optimal solution of the problem instance, as expected, while for the QUBO-formulated problem, we tested the open-source Julia libraries  {ToQUBO.jl} \cite{toqubo:2022}, Qiskit's optimization module \cite{quadratictoqubo}, the open-source Python library {PyQubo} \cite{zaman2021pyqubo} (in both cases operating under simulated annealing) and the emulation of the proprietary {Digital Annealer} of Fujitsu \cite{fujitsu}. 


\subsection{Quadratic Unconstrained Binary optimization}

The QUBO model can be applied to a wide range of combinatorial optimization problems that are known to be \NP-hard, such as the maximum cut, minimum vertex cover, multiple knapsack, and graph coloring problems. Its applications span a diverse set of domains, including the automotive industry \cite{glos2022optimizing}, portfolio optimization \cite{palmer2021quantum,9369145,Braine2021}, traffic flow optimization \cite{Braine2021}, job scheduling \cite{venturelli2015quantum,Zhang2022,Amaro2022}, railway conflict management \cite{domino2023quantum}, bioinformatics \cite{Matsumoto2022}, and others \cite{luckow2021quantum}. An extensive list of QUBO formulations for interesting problems can be found in \cite{ratke}.

Due to its one-to-one mapping to the Ising Hamiltonian, QUBOs have become a fundamental element of quantum-inspired computing. Both the Digital Annealer developed by Fujitsu and the adiabatic quantum computers manufactured by D-Wave Systems (as well as other vendors such as Qilimanjaro) employ the QUBO model to address complex optimization problems. Additionally, several approaches based on tensor networks seem to be suitable for a variety of QUBO-formulated optimization problems \cite{evenbly2011tensor,bridgeman2017hand,Zheltkov2020,nikitin2022quantum}. Although QUBO is particularly well-suited to these technologies, it can also be employed in NISQ devices using algorithms such as the QAOA. As such, the QUBO model is an important tool for quantum optimization with potential applications across a range of quantum computing platforms, and formulating constrained problems as such highly affects the quality of the solutions obtained.

Let us summarise the basics of QUBO via a graph problem. Given an undirected graph $G = (V, E)$  with a vertex set $V = \{1,2,...,N\}$ connected by the edge set $E=\{(i,j), \}$, $i,j\in V$, the cost function is defined as:
\begin{equation} \label{QUBO}
    \begin{split}
    & \min \quad \sum_{i=1}^N A_{ii}x_i + \sum_{i=1}^{N-1}\sum_{j>i}^{N} A_{ij}x_ix_j,
    \end{split}
\end{equation}    
where $x \in \{0,1\}$ are the binary variables and the elements $A_{ij} \in \mathbb{R}^{N \times N}$ are the problem instance parameters.

 At its most fundamental level, a QUBO can be expressed as:
\begin{equation}
    \begin{aligned}
         \min & \quad x^T Qx + b,
    \end{aligned}
\end{equation}
{where the decision matrix $Q \in \mathbb{R}^{N \times N}$ contains the problem instance and $b\in \mathbb{R}$   is a constant offset term. \par}

By using a suitable change of variables $x_i = \frac{1 - \sigma_i }{2}$,  Eq.(\ref{QUBO}) can be mapped onto the Ising Model Hamiltonian as:

\begin{equation}
    H = - \sum_j h_j\sigma_j - \sum_{j<k} J_{jk}\sigma_j\sigma_k,  
\end{equation}
where $\sigma\in  \{-1,1\}^N$ are the (classical) spins, $h\in \mathbb{R}^N$ is the magnetic field, and $J \in \mathbb{R}^{N\times N}$, ${\rm diag}(J) = 0$, the spin-spin interaction symmetric matrix between adjacent spins $j$ and $k$. See App. \ref{ap:qubo_to_ising} for more details. The problem to be solved then is:
\begin{align}
    \min_{\sigma_i \in \{-1,1\}} & \quad H 
\end{align}

For the QUBO formulation of \KP~we can take several slack-based approaches including ``off-the-shelf'' LP-to-QUBO converters such as Qiskit's  {QuadraticProgramToQubo} class and methods as well as the Julia package  {ToQUBO.jl}. Further to that, we can perform a ``custom'' slack formulation and also use the unbalanced penalization method we mentioned previously. 

\subsection{Slack Variable Formulation}

In the process of converting MILPs to QUBOs, it is common practice to introduce a slack variable, $S \in \mathbb{R}_{\geq 0}$ (whose purpose will be discussed shortly), for each linear inequality and transform it into an equivalent linear equality. Subsequently, a penalty term is constructed based on the slack variable, and the term is squared as per the standard approach, as outlined in \cite{glover2018tutorial} (see also \cite{symons2023practitioners}).

A variety of different slack-based QUBO formulations exist for the \KP~\cite{quintero}. Here, the corresponding penalisation term with weight $\lambda_{0}\in \mathbb{R}_{+}$ is given by the equality:
\begin{align}
    \lambda_{0}\left(\sum_{i=1}^{n} w_{i}x_{i} - W + S\right)^{2} = 0.
\end{align}
The purpose of the auxiliary slack variable $S$ is to reduce this term to $0$ once the constraint has been satisfied, $0\leq S \leq \max_{x} \sum^{n}_{i} w_{i}x_{i} - W$. In practice, $S$ is decomposed into binary representation using variables $s_{k}\in \{0,1\}$ as follows:
\begin{align}\label{eq:slack}
    S = \sum^{N_{s}}_{k=1} 2^{k-1}s_{k}.
\end{align}
The parameter $N_{s}$ corresponds to the number of binary variables required to represent the maximum value that can be assigned to the slack variable,
and in the case of \KP, $N_{s} =  \lceil \log_{2}(W)  \rceil$, where $\lceil x  \rceil$ is the ceiling function. Formulating the slack variable as in Eq. \eqref{eq:slack} is commonly referred to as the ``log-encoding''. 

The full QUBO formulation for the \KP~takes the form of maximizing the objective function:

\begin{align}\label{eq:KPslack1}
     \sum^{n}_{i}v_{i}x_{i} - \lambda_{0}\left(\sum_{i}^{n} w_{i}x_{i} - W + \sum^{N}_{k=1} 2^{k-1}s_{k}\right)^{2},
\end{align}
which can be understood as an augmented Lagrangian \cite{Bertsekas1982}. Alternatively, there are other QUBO formulations of the \KP~that we could decide to use that follow a similar Lagrangian paradigm. For example, we can instead consider maximizing the following objective function:
\begin{align}
   \sum^{n}_{i}v_{i}x_{i} - \lambda_{0}\Big(\sum_{i}^{n} w_{i}x_{i} -  \sum^{W}_{k=1} ks_{k}\Big)^{2} - \lambda_1 \Big(1 - \sum^{W}_{k=1}s_{k}\Big)^{2}.
    \label{eq:KPslack2}
\end{align}
In this formulation, known as ``one-hot encoding", the number of slack bits is equal to the capacity of the knapsack, $W$. Here, an additional penalty term is required to enforce only one of these slack bits to be assigned a value of 1. A drawback of this formulation is that the binary input length for the slack variables scales linearly with the values of the constraints, hence it can lead to an unreasonably large number of bits for problem instances with large $W$ which can exhaust available resources. This issue becomes quite relevant in Sec. \ref{sec:illustrations}.

\subsection{Balanced Slack-Based Approaches}\label{sec:encodings}
In this section, we give an overview of the approaches used to determine solutions for different balanced formulations of the \KP~QUBO. The known optimal solution for our small instance that we consider corresponds to an objective value of $309$ and uses the full capacity of the knapsack. 


\vspace{1em}

\noindent
\emph{\textbf{Off-The-Shelf Converters}}. The first ``off-the-shelf'' approach we tried is {ToQUBO.jl}, an open-source Julia package that automatically reformulates a variety of optimization problems, including MI(L)Ps, to a QUBO. The user can use the  {JuMP} \cite{DunningHuchetteLubin2017} package to build the MILP form of \KP.  {ToQUBO.jl} provides 6 ways for encoding variables into binary representations. This includes the logarithmic and one-hot approaches mentioned above as well as other less well-known techniques. For continuous decision variables, this can be very useful since the user can provide a tolerance factor to manage the upper bound on the representation error caused by the binarization. Additionally, {ToQUBO.jl} works in conjunction with {QUBODrives.jl}, a companion package that provides common API to use QUBO sampling and annealing machines such as D-Wave's simulated annealer and, with license, the quantum annealer via {DWaveNeal.jl}. We make use of {ToQUBO.jl} to employ both of the aforementioned binary encodings, both successfully finding the optimal solution. However, the unbalanced QUBO formulation is a recent development which has not been widely adopted, hence its encoding is not available as part of the {ToQUBO.jl} package. 

Another ``off-the-shelf'' converter is provided by Qiskit's optimization module which includes functionality for automatically transforming quadratic programs into QUBOs (the binary property allows us to use such). This transformation can be accomplished by first initializing a  {QuadraticProgram} and subsequently utilizing the  {QuadraticProgramToQubo} class to convert it into a QUBO via the log-encoding method for slack variables. The module allows the formulated QUBO to feed into several algorithms used by Qiskit to solve optimization problems, such as {SamplingVQE} or {QAOA}, however, the user can also extract the coefficient matrix to use with other solvers. We input this coefficient matrix into {neal} and Fujitsu's Digital Annealer, where we found that both solvers are able to reach optimum, with the caveat that a large number of runs are needed compared to other methods. For larger problem instances, this can become computationally expensive and thus may be an inadequate choice for \CO.

There are several variations of the balanced slack-based approach to QUBO formulations, which we summarize now.
\vspace{1em}

\noindent
\emph{\textbf{Log-encoding}}. This approach refers to implementing Eq. \eqref{eq:KPslack1}. 
The number of slack bits required for this approach, for the instance of interest, is $N_{s} = 8$. Two different regimes for the weight of the penalty term, $\lambda_{0}$, are checked: 
\begin{enumerate}[label=(\alph*),font=\itshape]
    \item where the penalty term and the cost function have equal weighting ($\lambda_0 = 1$), and
    \item where the penalty term is more important than the cost function ($\lambda_0 = 1\times 10^{4}$).
\end{enumerate}
We utilized the {neal} package to implement D-Wave's simulated annealer as our heuristic optimizer in both scenarios. The simulated annealer successfully returned the optimal solution in both cases, which was consistent with the results obtained from the classical solvers. The (emulation of the) Digital Annealer from Fujitsu was employed for both scenarios as well. Similar to the case with {neal}, it successfully produced optimal solutions in both instances.





\noindent
\emph{\textbf{One-hot encoding}}. 
We repeated the previously described process but with a focus on analyzing the solutions of Eq. \eqref{eq:KPslack2}. As previously mentioned, this type of formulation requires a large number of bits to encode the slack variables for large knapsack capacities, in this case, 165 bits. In this approach, we once more explored various weight regimes for $\lambda_0$ and $\lambda_1$, in relation to the weight of the cost function term. By selecting $\lambda_0 = 10^{-1}$ and $\lambda_{1} = 10^{3}$, we consistently identified the optimal solution for both {neal} and Fujitsu cases.



\subsection{Unbalanced Penalization Approach}\label{sec:unbalanced}

Given that the number of qubits required scales proportionally with the number of variables, we sought to employ the methodology proposed by Montanez-Barrera \emph{et. al.} \cite{montanez2022unbalanced} in order to eliminate the need for slack variables. To accomplish this, we adopted an approximation technique that creates penalty terms that take on small values when the constraint is fulfilled and large values when it is violated. We rearrange the inequality to define an auxiliary function:
\begin{align}
    h(x) = \sum^{n}_{i}w_{i}x_{i} - W \leq 0.
\end{align}
{ Using the exponential function $f(x) \coloneqq e^{h(x)}$. 
the constraint on $h(x)$ is satisfied. 
However, since only linear and quadratic terms may be encoded into a QUBO}, it is necessary to use a 2\textsuperscript{nd} order Taylor approximation of $f(x)$. For weights $\lambda_{0}$ and $\lambda_{1}$, the resulting QUBO for the \KP~problem reads:

\begin{align}
     \sum^{n}_{i}v_{i}x_{i} + \lambda_{0}\Big(\sum_{i}^{n} w_{i}x_{i} - W\Big) + \lambda_{1}\Big(\sum_{i}^{n} w_{i}x_{i} - W\Big)^{2}  = 0.
\end{align}

It should be noted that Ref. \cite{montanez2022unbalanced} investigated this type of QUBO formulation for different instances of the Traveling Salesperson Problem (\TSP), Bin Packing Problem (\BinPP), and \KP~and found that the minimum energy eigenvalue of the corresponding Hamiltonian did not necessarily coincide with the optimal solution. However, the optimal solution was always found to be amongst the lowest-energy eigenvalues which enhances the confidence for this approach in large-scale experiments.

To clarify, the advantage of reformulating the problem as an unbalanced QUBO rather than the traditional balanced approaches is that one significantly reduces the number of variables, and hence the bits required to represent the problem, which in effect reduces the resource cost as well as the search space of the optimal solution. However, drawbacks associated with the unbalanced approach are that, because of the use of a heuristic penalisation function our constraints are less strict and the groundstate of the corresponding Ising Hamiltonian is less likely to coincide with the optimal solution of the problem at hand.

{ Using the weights provided by Ref. \cite{montanez2022unbalanced}, as a quick check, we were able to reach  optimality for this \KP~instance for both PyQUBO and Fujitsu unbalanced approaches, signifying that the unbalanced formulation exhibits (some) robustness. For larger \KP~instances tested we produced close to optimal results which, however, would periodically softly break the maximum weight limit.}



\subsection{Summary of the simulations}

Below, in Table \ref{table:1}, we provide a list of MILP and QUBO solvers, using the two aforementioned approaches, slack-based and unbalanced.

\begin{table}[!htb] 
\begin{tabular}{ccc}\toprule[2pt]
Problem Encoding & Reference & Solver    \\ \hline 
\multirow{2}{*}{ILP} &\cite{GLPK}& GLPK \\
                     &\cite{Huangfu2017}& HiGHS \\ \Xhline{0.02\arrayrulewidth}
\multirow{1}{*}{ {ToQUBO.jl}} &{\cite{toqubo:2022}}&  {D-Wave} \\
                           \Xhline{0.02\arrayrulewidth}
\multirow{2}{*}{ {QuadraticProgramToQUBO}} &\multirow{2}{*}{\cite{quadratictoqubo}}& D-Wave (PyQubo)\\
                           && Fujitsu \\ \Xhline{0.02\arrayrulewidth}
\multirow{2}{*}{Log} &\multirow{2}{*}{Sec. \ref{sec:encodings}} & D-Wave (PyQubo) \\
                           && Fujitsu\\ \Xhline{0.02\arrayrulewidth}
\multirow{2}{*}{One-hot} &\multirow{2}{*}{Sec. \ref{sec:encodings}}&   D-Wave (PyQubo)  \\
                           && Fujitsu \\ \Xhline{0.02\arrayrulewidth}
\multirow{2}{*}{Unbalanced} &\multirow{2}{*}{Ref. \ref{sec:unbalanced}, \cite{montanez2022unbalanced}}  & D-Wave (PyQubo) \\
                          & & Fujitsu  \\ \Xhline{0.02\arrayrulewidth}
\end{tabular} 
\caption{Summary of methods used to solve the provided instance of the \KP. All methods found the optimal solution.
}
\label{table:1}
\end{table}



\noindent

\subsection{Survey of Alternative Approaches}

There exist various alternative approaches that exhibit varying degrees of divergence, in terms of algorithmic implementation, from the aforementioned methods. This section serves solely to provide a survey of some of these approaches without undertaking an experimental analysis. 
\vspace{0.2em}

\noindent
\emph{\textbf{Quantum Hybrid Frank-Wolfe method}}.
Recently, a hybrid quantum generalization to the Frank-Wolfe method was proposed in \cite{yurtsever2022q}. This hybrid (Quantum) Frank Wolfe (Q-FW) augmented Lagrangian
method is suitable to tackle large QUBO instances due to a tight copositive relaxation of the original QUBO formulation while dealing with the expensive hyper-parameter tuning found in other QUBO heuristics. The Q-FW method first formulates constrained QUBOs as copositive programs, then employs the Frank-Wolfe method while satisfying linear (in)equality constraints. This is the converted to a set of unconstrained QUBOs suitable to be run on, e.g., quantum annealers. 
It was found that Q-FW
successfully satisfied linear equality and inequality constraints, in the context of QUBOs in computer vision applications and the authors of \cite{yurtsever2022q} solved intermediary QUBO problems on actual quantum
devices demonstrating that Q-FW offers a promising alternative the validity of Q-FW, to ``traditional'' quantum QUBO solvers.

In a broader context, Q-FW seems to have the ability to address the costly hyper-parameter tuning associated with other QUBO heuristics. By formulating constrained QUBOs as copositive programs, it adeptly handles linear equality and inequality constraints and transforms them into unconstrained QUBOs compatible with quantum annealers or other Ising machines. The general applicability and comparative efficiency of Q-FW against established quantum QUBO solvers is still unclear and subject to further research.

\vspace{0.2em}

\noindent
\emph{\textbf{Grover adapted binary optimization}}.
Moving to Fault-Tolerant architectures, quadratic speed-up for combinatorial optimization problems is achievable with the Grover adapted binary optimization (GABO) \cite{gilliam2021grover} when compared to brute force search. However, to achieve this, efficient oracles must be developed to represent problems and identify states that satisfy specific search criteria. Quantum arithmetic is commonly utilized to accomplish this task, but this approach can be expensive in terms of required Toffoli gates and ancilla qubits, which may pose a challenge in the near future. Interestingly \cite{gilliam2021grover} provides such an oracle construction which makes GABO a promising approach for future quantum computers. 

GABO might be able to offer a significant quantum advantage, a Grover-like quadratic speed-up for combinatorial optimization problems compared to traditional brute force search. However, this speed-up can be only impactful in the realm of Fault-Tolerant quantum architectures, and thus not applicable to NISQ devices.

\vspace{0.2em}

\noindent
\emph{\textbf{Graph Neural Networks}}.
Physics-inspired Graph Neural Networks (GNNs) have been used in \cite{Schuetz2021} (see \cite{velivckovic2023everything} for a survey) where a physics-informed GNN-based scalable general-purpose QUBO solver is proposed. The approach therein is suitable for encoding any $k$-local Ising model such as the $k=2$ \CO~ problem discussed later. The GNN solver first drops the integrality constraints in order to obtain a differentiable relaxation $f'$ of the original objective function $f$ and subsequently proceeds to unsupervised learning on the node representations. The GNN is then trained to generate soft assignments, predicting the likelihood of each vertex in the graph belonging to one of two distinct classes, in conjunction with heuristics that aid in the consistency of the problem. Interestingly, the authors benchmark this approach in demanding problem instances of \MAXCUT~to find that it competes with the best in class SDP algorithms such as the Goemans-Williamson algorithm \cite{Goemans1995}. 

The parallel processing capabilities of modern GPUs are well-suited for GNN operations, making it feasible to handle large-scale graphs and complex optimization problems. Furthermore, the adaptability of GNNs, combined with unsupervised learning, allows for general-purpose solutions that can be applied across a variety of problem instances without the need for extensive retraining. This universality potentially saves computational resources in the long run. However, possible relaxations of integrality constraints for differentiability can sometimes lead to solutions that are not directly applicable or optimal for the original discrete problem. Another remark is that while GNNs can predict soft assignments efficiently, converting these into definitive solutions might lead to suboptimal results.



\vspace{0.2em}

\noindent
\emph{\textbf{QUBO continuous relaxations with light sources}}. Finally, let us mention a recent heuristic quantum-inspired (relaxation) approach for solving QUBOs, introduced in \cite{meirzada2022lightsolver}. Concretely, the binary variables of the QUBO problem are represented by the relative phases of laser sources, transforming the discrete optimization problem into a continuous one. The lasers interact through a unique optical coupler, which uses programmable diffractive elements and additional optical components to control the interaction between all pairs of lasers, with a dynamic range of up to 8 bits. This design enables a fully connected network between all lasers, facilitating high-resolution pairwise interactions that are crucial in solving QUBO problems.

In \cite{meirzada2022lightsolver} a benchmarking was performed for instances of the 3-regular 3-\XORSAT problem and it was found that this method achieves significantly better TTS results as the instance size was increased. Despite this, seemingly incredible result, it is quite unclear whether this approach can perform on par for problems that do not have known poly-time solutions (for generic instances of 3-\XORSAT~there exists an algorithm with $\mathcal{O}(N^{2.736\ldots})$ complexity). In fact, it would be interesting to benchmark this approach against specialized SAT and MIP solvers. 

 Among its advantages, the method's transformation of discrete variables into continuous may offer new avenues for efficient problem-solving. The fully connected network ensures precise pairwise interactions, a fundamental requirement for many combinatorial problems, also one that not all Ising machines can achieve. The aforementioned results of this method are obtained via an emulation architecture, rather than the physical machine which in turn might be able to perform even better. Another benefit of this tech is that it is readily available. However, there are a number of uncertainties. While the method showed promise on specific benchmarks like the 3-\XORSAT problem, its general applicability and competitiveness against established solvers remain untested and benchmarking e.g. within the MaxSAT Evaluation challenge would be highly en lightning. Furthermore, its performance on problems lacking known polynomial-time solutions is yet to be determined, raising questions about its versatility.

\vspace{0.2em}

\noindent
\emph{\textbf{Simulated quantum annealing}}.
Obtaining quantum advantage for certain  problem instances may be, in some cases, closer than one would anticipate by utilizing the technique of simulated quantum annealing (SQA) \cite{crosson2016simulated}. Specifically, in \cite{crosson2016simulated} a method to rigorously demonstrate that the Markov chain underlying SQA effectively samples the target distribution and discovers the \emph{global minimum} of the spike cost function in \poly~time supporting is developed. While the analysis is limited to a very specific model and cannot be considered conclusive, the authors use interesting techniques such as initiating warm starts (a very popular technique applied in Deep Learning \cite{loshchilov2016sgdr} as well as in QAOA \cite{egger2021warm}), from the adiabatic path and using the quantum ground state probability distribution to comprehend the stationary distribution of SQA. 

Interestingly, SQA might be considered effective, as compared to a purely classical algorithm, for optimization problems with spike (deep and narrow) global optima. As such, exploiting hybrid solvers that combine SQA and classical algorithms could offer an alternative approach \cite{10078403} due to SQA's effectiveness in sampling target distributions. However, to the best of our knowledge, its applicability is currently limited to specific models, and the broader effectiveness remains an open question.


\section{Collateral optimization }
\label{sec:level2}

In Section \ref{sec:knapsack}, our main objective was to identify the most suitable formulations of the \KP~as a QUBO. However, for smaller-scale problems, MILP solvers are generally expected to perform better than heuristic, hybrid, and near-term quantum solvers. Therefore, in this section, we first formulate the Collateral Optimization problem \CO~as a MILP and subsequently re-formulate it as a QUBO, similar to the approach taken in Section \ref{sec:knapsack}.

The objective is to conduct several small-scale problems utilizing hybrid solvers with the purpose of algorithmic reformulation of the problem and its implementation. The ultimate goal is not to demonstrate the potential of hybrid or quantum alternatives for combinatorial optimization, but rather to establish an automated approach to solving the problem once hardware capabilities become more advanced. In Table \ref{tab:TableofDefinitions} we describe the most common financial terms used in the context of \CO, while Table \ref{tab:TableOfNotation} summarises the mathematical notation that will be utilized in the following.


\subsection{\CO~MILP Formulation}

In order to mitigate the risk of a borrower defaulting on a loan, it is necessary for them to furnish collateral in the form of stocks, bonds, cash, or other assets to offset any outstanding exposure. In the present scenario, we consider a financial institution that has a collection (or inventory) of assets, indicated by $\mathcal{I}$, which must be allocated among a set of accounts, indicated by $\mathcal{A}$. We shall use the indices $i$ and $j$ to refer to an asset and an account, respectively. The total number of assets and accounts is given by $n$ and $m$, respectively. It should be noted that our primary focus in this section is on the algorithmic reformulation of the problem and its implementation, which could be automated once sufficiently powerful hardware becomes available. Consequently, we shall run a set of small-scale problems employing hybrid solver emulators as in Sec. \ref{sec:knapsack}. Our objective is not to demonstrate the superior performance of hybrid or quantum alternatives for combinatorial optimization but rather to heuristically identify the most appropriate formulation of the problem at hand.

\begin{figure}[htb!]
    \centering
    \includegraphics[scale=0.14]{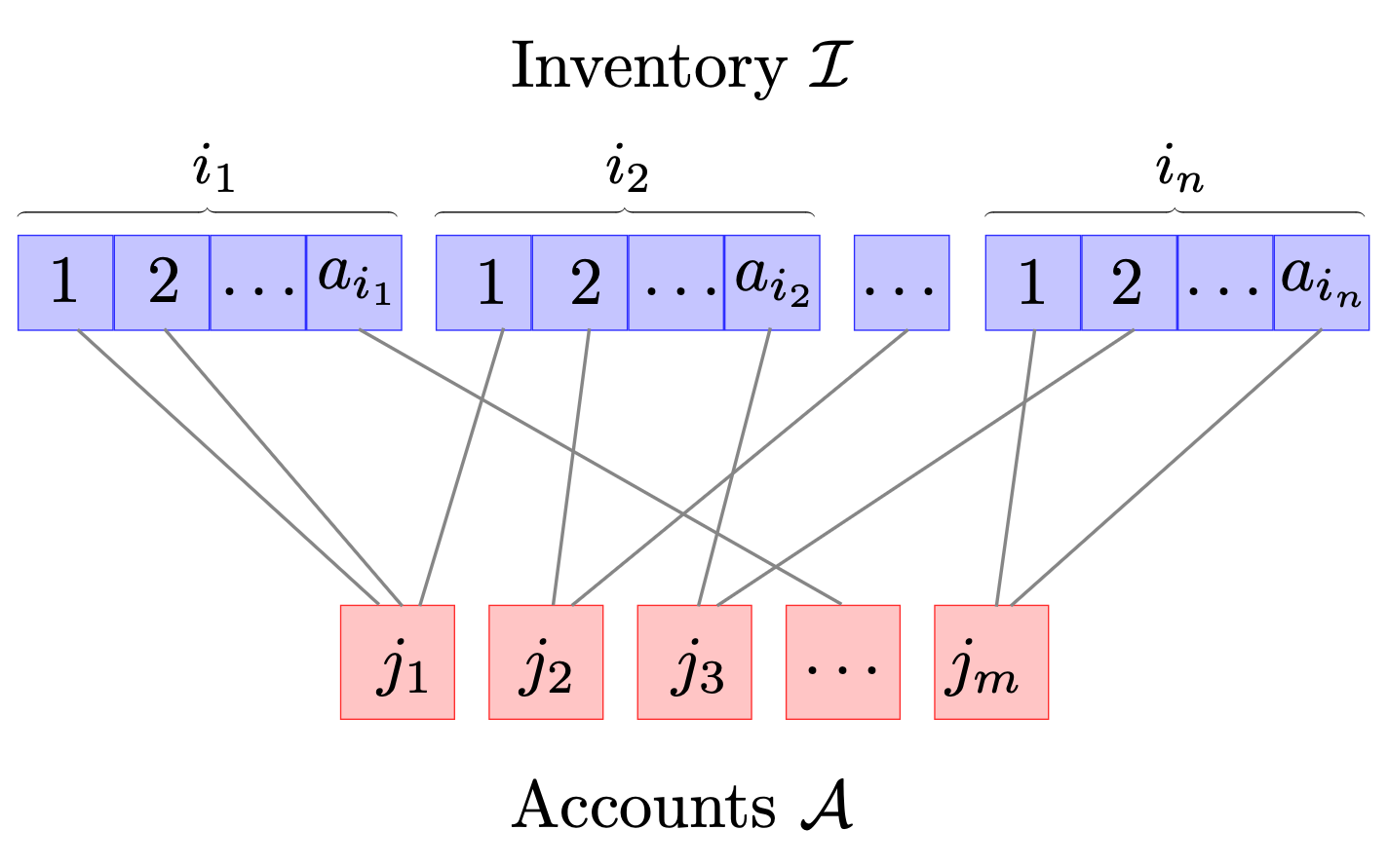}
    \caption{Schematic representation of the collateral optimization problem. The figure above can be seen as a bipartite graph where one has to make optimal allocation of non-unique and weighted pairings.}
    \label{fig:fig1}
\end{figure}

Interestingly, \CO~can be formulated as a bipartite matching problem as in Fig. \ref{fig:fig1}. This bipartite graph is created with two sets of nodes: one set representing the inventory of assets $\mathcal{I}$, and the other set representing accounts $\mathcal{A}$. The edges between these nodes represent potential allocations of assets to accounts, with weights on these edges representing the suitability, cost, or value of the allocation, and the edges are multi-directed and need to respect certain constraints. To model the constraints of \CO~modifies the graph accordingly, as we will explain below.

\emph{\textbf{Limits}}. Delving deeper into the problem, consider each asset $i_k \in \mathcal{I}$, where $\mathcal{I} = {i_1, i_2, \ldots, i_n}$. Each asset $i$ (momentarily simplifying the indices) is subdivided into a maximum quantity denoted by $a_{i}$. In more formal terms, there is a constraint on the maximum quantity of asset $i$ that can be assigned. This is useful because, in the context of utilizing stocks as collateral, a financial institution may need or require the enforcement of an upper limit regarding the number of shares that can be allocated. The quantity of asset $i_k$ can be converted into a corresponding dollar value by multiplying by the dimensionful term $v_{i}$, which is the market value (USD) per unit quantity. 


\emph{\textbf{Tiers}}. { Every asset is linked to a tier, represented as $\omega_{i} \in [0,1]$}. This tier acts as a measure of the asset's quality, where distinct tiers correspond to various degrees of quality or attractiveness in the context of the \CO~problem. 
The higher the value of $\omega_{i}$, the higher the quality of asset $i$.

\emph{\textbf{Exposure}}. When a financial institution borrows from one of its lenders, collateral must be posted to adequately cover capital that could be lost in the event of a default. This capital requirement is known as ``exposure''. For each account $j$, there is a required exposure (in USD) that must be met indicated by $c_{j}$. Additionally, the duration of a transaction to a particular account can either be short term or long term. A binary variable $d_{j}\in \{0,1\}$ is used to indicate the duration of account $j$. A value $1$ is assigned to short term and $0$ for long term transfers. To reduce the risk of losing posted collateral, it is required to minimize the use of high quality assets for long term transactions whilst maximizing their use for short term transactions. 
We chose our decision variable to be a matrix $Q\in \mathbb{Q}_{ \leq 1}^{n\times m}$, where the element $Q_{ij}$ is the fractional amount of asset $i$ that is allocated to account $j$. That is, the rows of $Q$ correspond to the number of assets and the columns of the partition of each. 
\begin{equation}
    Q = \begin{pmatrix}
            a_{11}  & \ldots & a_{1m} \\
              \vdots   & \ddots &  \vdots   \\
            a_{1n}  & \ldots & a_{nm} \
        \end{pmatrix}.
\end{equation}

To illustrate an objective function that represents our goals consider the simple case of allocating a single asset, $i$, across two accounts, $j$ and $l$, which have long and short term requirements, respectively. In this case, the objective function can be formulated as:

\begin{equation}
    \min_{} \quad \omega_{i} Q_{ij} + (1 - \omega_{i})Q_{il}.
    \label{eq: simpleObjectiveFunction}
\end{equation}
In the expression above, the coefficient preceding short term allocations is set to $1 - \omega_{i}$. This is so that we favor allocations of higher quality assets for trades with a short duration. To generalize this for all assets and accounts, we need a mechanism which updates these tiers according to the type of account they are posting collateral towards. This can be done by constructing a coefficients matrix $\Omega$, where each element can be determined using:

\begin{equation}
\Omega_{ij} =  |\omega_{i} - d_{j}|.
\end{equation}

Just as with the tiers, each element of $\Omega$ has a range of $[0,1]$. The objective function is then just the sum of element-wise multiplications between $\Omega$ and $Q$:

\begin{equation}
    \min_{Q} \quad \sum_{i=1}^{n}\sum_{j=1}^{m} \Omega_{ij} Q_{ij}.
    \label{eq:objectiveFunction}
\end{equation}


To post collateral such that the financial institution meets the exposure for each account, we include a requirement constraint:

\begin{equation}
    \sum_{i=1}^{n} Q_{ij}a_{i}v_{i} H_{ij} \geq c_{j}, \quad \forall j \in \mathcal{A}.
    \label{eq:requirementConstraint}
\end{equation}
Here, $v_{i}$ denotes the dollar market value for a single unit of quantity for asset $i$. Hence,
the term on the left-hand side represents the dollar value for the quantity of collateral that is chosen to be allocated, adjusted by a fractional factor $H_{ij}$ which is referred to as the haircut. Since markets are dynamic, the value of a posted collateral can diverge from its market value over time. In the case that the value drops below the required collateral value, the receiver is at risk. To avoid this, each account owner can evaluate the risk and place a haircut factor to reduce the value of an asset. The haircut is defined as the percentage difference between the market value and its value whilst used as collateral. For example, a haircut of 10\% corresponds to  $H_{ij} =0.9$, meaning the collateral value is 90\% of the original market value. 

We need to ensure that we do not allocate more collateral than we have available in inventory (i.e. we do not allocate more than 100\% of the maximum available quantity). In financial terms, this prevents us from short selling the asset. This is done by including a consistency constraint:

\begin{equation}
    \sum_{j=1}^{m} Q_{ij} \leq 1 \quad \forall i \in \mathcal{I}.
    \label{eq:consistencyConstraint}
\end{equation}

There is also the trivial constraint to make sure that $Q_{ij}$ does not take negative values:

\begin{equation}
    Q_{ij} \geq 0 \quad \forall i \in \mathcal{I}, \forall j \in \mathcal{A}.
    \label{eq:nonnegativeConstraint}
\end{equation}

\noindent
\emph{\textbf{Further constraints that be can be imposed.}} For example, there may be limits on the amount of a particular asset $i$ to account $j$, given by $B_{ij}$. If $B_{ij} = 0$, the allocation is not eligible. This is a one-to-one constraint and has the following mathematical form:
\begin{equation}
    Q_{ij}a_{i} \leq B_{ij} \quad \forall i \in \mathcal{I}, \forall j \in \mathcal{A}.
    \label{eq:oneToOneConstraint}
\end{equation}


Moreover, there are constraints that restrict the allocation of specific groups of assets to a single account, which exhibits a many-to-one relationship. For instance, certain types of assets, say $\{i_{X_1},i_{X_2},i_{X_3}\}$, may be subject to restrictions due to their interrelationships (for example there exists a parent company $X$ that posts these assets). To formalize this constraint, we introduce $\mathcal{G}$, which represents the set of all groups of assets. The binary variable $T_{ig}$ is used to indicate whether asset $i$ belongs to the group $g$, while $K_{gj}$ represents the upper bound on the total amount of assets from group $g$ that can be allocated to account $j$:


\begin{equation}
    \sum_{i=1}^{n} T_{ig}Q_{ij}a_{i} \leq K_{gj}, \quad \forall g \in \mathcal{G}, \forall j \in \mathcal{A}.
    \label{eq:manyToOneConstraint}
\end{equation}


The aim of imposing limits on the allocation of assets is to promote diversification and thereby reduce the risk borne by the receiver. In this paper, we concentrate on allocating cash rather than equity and bonds, which allows us to avoid the constraint that $Q_{ij}a_{i}$ must take an integer value. Consequently, we can formulate the problem of collateral optimization as a continuous optimization problem without the need for additional constraints.

\vspace{1em}

\noindent
\emph{\textbf{The Complete \CO~ Problem}}. Taking into account the information presented in the previous paragraphs, we can express the \CO~problem using a MILP formulation:
\begin{subequations}
\begin{align} \label{eq:obj}
     {\min_Q} & \quad \sum_{i=1}^{n}\sum_{j=1}^{m}\Omega_{ij} Q_{ij}&\quad \forall i \in \mathcal{I}, \forall{j} \in \mathcal{A}  \\ \label{eq:con1}
    {\text{s.t.}}  & \quad \sum_{j=1}^{m} Q_{ij} \leq 1 & \quad \forall  i \in \mathcal{I}\\ \label{eq:con2}
    & \quad  Q_{ij}a_{i} \leq B_{ij} &\quad \forall i \in \mathcal{I}, \forall j \in \mathcal{A} \\ \label{eq:con3}
    & \quad \sum_{i = 1}^{n} T_{ig}Q_{ij}a_{i} \leq K_{gj} &\quad \forall g \in \mathcal{G}, \forall j \in \mathcal{A} \\ \label{eq:con4}
    & \quad  \sum_{i = 1 }^{n} Q_{ij}a_{i}v_{i}H_{ij} \geq c_{j} &\quad \forall j \in \mathcal{A} \\ \label{eq:con5}
    & \quad Q_{ij} \geq 0 &\quad \forall i \in \mathcal{I}, \, \forall j \in \mathcal{A}.
\end{align}
\end{subequations}
For clarity, we will summarize the constraints presented above. Constraint \eqref{eq:con1} ensures that no asset is distributed to the accounts beyond unity. Constraint \eqref{eq:con2} amounts to the limit constraints for each asset-account pairing. Constraint \eqref{eq:con3} limits the quantity of particular groups of assets to certain accounts but will be ignored in what follows. Constraint \eqref{eq:con4} is the requirement constraint that enforces that we allocate a suitable value such that the lender's loan is secured. 
\\

\subsection{\CO~ QUBO Formulation}


\noindent
\emph{\textbf{Binarization}}. To formulate the QUBO, we need to make a change of variables so that our decision variable is represented by binaries. This change imposes certain limitations on the allocation of assets, which are discussed in detail below.
 We make use of an $n$-bit binary variable, similar to the methods used by Lang \cite{lang2022strategic} and Ottaviani \cite{ottaviani2018low}. To enable binary encoding of the decision variable $Q$, we can represent it as a matrix $q$ containing binary elements. This transformation enables us to only allocate assets in a limited number of ways, as detailed below: 
\begin{equation}
    q = \begin{pmatrix}
    q_{11} & \ldots & q_{01} \\
    \vdots & \ddots & \vdots \\
    q_{n1} & \ldots & q_{nm}
    \end{pmatrix}.
\end{equation}
Here, $q_{ij}$ is an $n$-bit binary variable expressing the fractional allocation of asset $i$ to account $j$:
\begin{equation}
    q_{ij} = x^{b=1}_{ij}, \ldots, x^{b=B}_{ij}, \quad x^b_{ij} \in \{0,1\},
    \label{eq:bitstring}
\end{equation}
where $B$ is the number of bits chosen. \par
Using $B = 4$ as an example, the largest number which can be represented by four bits is $1111_{\bin} = 15_{\dec}$. Thus, we split our allocation into 15 fractions, where if $q_{ij} = 0100_{\bin} = 4_{\dec}$ then 4/15 of asset $i$ is allocated to account $j$. By increasing $B$, we increase the precision of our allocations.

\par\par
By implementing a similar method to Braun \emph{et. al.} \cite{braun2023towards}, we can discretize our fractional allocation by discretizing the interval $[Q^{\min}_{ij}, Q^{\max}_{ij}]$:
\begin{equation}
   Q^{\max}_{ij} - Q^{\min}_{ij} = \sum_{b=1}^B p_{ijb} = \sum_{b=1}^B \frac{2^{(b-1)}(Q_{ij}^{\max} - Q_{ij}^{\min})}{M}. 
\end{equation}

$M$ is the maximum value which can be represented by a binary string of length $B$ ($M = 2^{B} - 1$). 
\emph{Remark}: $Q^{max}_{ij}, Q^{min}_{ij}$ have been included here for generality, but in our problem $Q_{ij}^{\max} =1 $ and $Q_{ij}^{\min} = 0$ by definition. Thus, $Q_{ij}^{\max} - Q_{ij}^{\min}$ is always $1$.

The discretized amount of item $i$ that is allocated to account $j$ is then described by the dot product of the binary vector $x_{ijb} = (x_{ij1}, \ldots, x_{ijB})^{T}$ and $p_{ijb} = (p_{ij1}, \ldots, p_{ijB})^{T}$. 
\\

\begin{align}
    q_{ij} = \sum^{B}_{b=1} p_{ijb}x_{ijb} \equiv p_{ijb}x_{ij}^{b},
\label{eq:binarized_matrix_elements}
\end{align}
where we use the Einstein summation convention (upper-lower repeated indices contract) for brevity. By pairing each bit in the bit-string of $q_{ij}$ with a coefficient $p_{ijb}$ allows us to represent more values in the allowed range which improves accuracy.

The cost function is then written as:
\begin{equation}
    \sum_{i=1}^{n}\sum_{j=1}^{m} 
    {\Omega}_{ij}p_{ijb}x_{ij}^{b}. \\
\end{equation}

We then make the same replacement of $Q_{ij} \rightarrow p_{ijb}x_{ij}^{b}$ in the remaining constraints of the problem to construct the binarized collateral optimization problem. Again, a wide number of classical solvers (open-source and commercial) are available which can find the global minimum of the problem in this form. 
\vspace{1em}

{ The total number of variables used to construct this binarized version of the problem is $\mathcal{O}(n m B)$. Since we have 
replaced the continuous variables with discrete binaries, the accuracy of the solution is expected to be reduced. This can be mitigated by increasing $B$, however, a compromise is needed between accuracy and resource usage. Regardless, the granularity for the fractional allocation is $1/M$.  }

\noindent
\emph{\textbf{Slack-based Formulation}}. 
We hereby introduce a QUBO formulation for the collateral optimization problem, incorporating slack variables. The constraints outlined in \eqref{eq:con1} - \eqref{eq:con4} significantly influence the number of bits necessary to encode these slack variables, potentially leading to an extensive bit requirement. To address this issue, we employ the log-encoding method earlier in Sec. \ref{sec:encodings}.

For constraints expressed as ``less-than-or-equal-to'' inequalities, the number of bits required to encode the slack variable can be readily computed using $\lceil \log_{2}u \rceil$, where $u$ represents the upper bound of the constraint. Nonetheless, addressing the requirement constraint \eqref{eq:con4} necessitates a more nuanced approach. As the objective is to minimize the excess value of the collateral posted, employing slack variables might prove inadequate, since their purpose is to diminish the corresponding penalty term to $0$ for any values satisfying this constraint. In contrast, within MILP frameworks, the solution to a minimization problem typically aligns closely with the lower bound of such a constraint. In light of this observation, we choose to alter the exposure requirement by transforming it into an equality constraint, thereby relaxing the original formulation:
\begin{align}
    \sum^{n}_{i=1} Q_{ij}a_{i}v_{i}H_{ij} = c_{j}, \quad \forall j \in \mathcal{A},
\end{align}
and as a result, the associated penalty term requires no slack variables.


Hence, with the objective of minimizing the associated penalty terms originating from the aforementioned expression, the QUBO should yield a solution conforming to the boundaries of the exposure constraints. A limitation of this strategy is the intrinsic stochastic characteristic of annealing techniques and their propensity to become ensnared in local minima, potentially resulting in marginally exceeding or not quite meeting the mandatory exposure. Furthermore, because the upper bound for each consistency constraint is equal to one, we need to adjust their form so that we can calculate the number of slack bits required. We proceed by rearranging the binarized version of this constraint which then attains a fractional form:
\begin{align}
    \sum^{m}_{j=1} p_{ijb}x_{ij}^{b}  = \sum^{m}_{j=1} \sum^{B}_{b=1} \frac{2^{b-1}  x_{ij}^{b}}{M} \leq 1,
\end{align}
$\forall i \in \mathcal{I}$. By multiplying both sides by $M$, it is straightforward to realize that the highest value that can be represented by the bitstring is the one we use as our new upper bound, and this further allows us to easily determine the number of slack bits needed. 

The penalty term for each of the $n$ consistency constraints is written as:
\begin{align}
\sum_{i=1}^{n}\left(\sum_{j=1}^{m}M(p_{ijb}x_{ij}^{b} -  1) 
 + S_{\text{con}} \right)^{2},   
\end{align}
where $S_{\rm con}$ is the slack variable for each constraint which is encoded by binary variables $s_{k}$ via:

\begin{align}
    S_{\text{con}} = \sum_{k=1}^{\lceil \log_{2}(M) \rceil} 2^{k-1}s_{k}.
\end{align}
Instead of introducing a penalty term for one-to-one constraints \eqref{eq:con2}, we can ensure that these are satisfied by reducing the number of bits representing each allocation so that these limits cannot be violated. The number of bits representing an allocation, $n_{ij}$, can be determined by:

\begin{align}
    n_{ij} = \left\lfloor \log_{2}\left(\frac{B_{ij}}{a_{ij}}M\right)  \right\rfloor.
\end{align}
Note that now in \eqref{eq:binarized_matrix_elements} the upper limit of the summation is no longer $B$ but now $n_{ij}$ as this is the number of bits representing the allocation $q_{ij}$.

We choose to floor the result from the logarithmic function, as the alternative would still allow violations. However, a consequence of this is that the one-to-one constraints in the QUBO are more restrictive than their MILP counterpart. 


Aside from these nuances above and the additional step of binarization, constructing the balanced formulation for this QUBO follows the same process described above in the discussion of the \KP~instance. For the many-to-one constraints \eqref{eq:con3}, we introduce log-encoded slack variables $S_{K_{ij}}$. 
It is straight-forward to derive the full QUBO objective as:

\begin{equation}
\begin{split}
   & \lambda_{0}\sum_{i=1}^{n}\sum_{j=1}^{m} 
    {\Omega}_{ij}p_{ijb}x_{ij}^{b}  \\
   & + \lambda_{1}\sum_{i=1}^{n}\Big(\sum_{j=1}^{m}M(p_{ijb}x_{ij}^{b} -  1) 
+ S_{\rm con} \Big)^{2} \\
    & + \lambda_{2}\sum_{j = 1}^{m}\Big(\sum_{i = 1}^{n}p_{ijb}x_{ij}^{b} a_{i}v_{i}H_{ij} - c_{j}\Big)^{2} \\
   & + \lambda_{3}\sum_{j=1}^{m}\sum_{g=1}^{G}\Big(\sum_{i=1}^{n}p_{ijb}x_{ij}^{b} T_{ig}a_{i} - K_{gj} + S_{K_{gj}}\Big)^{2}. \\
\end{split}
\end{equation}
\vspace{1em}

\noindent
\emph{\textbf{Unbalanced formulation}}. The previous constraints, set in equations \eqref{eq:con1} -\eqref{eq:con5}, can be converted into penalty terms for the QUBO through unbalanced penalization, as we displayed with the \KP~instance in Sec. \ref{sec:unbalanced}. For this approach, auxiliary functions are defined from the constraints, and Taylor approximations of appropriate exponentiations of these functions are used to derive the penalty terms.

For instance, consider the consistency constraint equation \eqref{eq:con1}. We move the upper bound to the LHS of the inequality and set this as our auxiliary function $h(x)$:
\begin{equation}
    h(x) = \sum_{j=1}^{m}p_{ijb}x_{ij}^{b} - 1 \leq 0.
\end{equation}
Since this is a `less-equal to zero' inequality we use $e^{h(x)}$ to derive a penalty term that takes small values when this constraint is satisfied and large values when violated. As mentioned previously, QUBOs may only contain linear and quadratic terms so, therefore, we need to take a 2\textsuperscript{nd} order Taylor approximation to obtain: 

\begin{equation}
    \lambda_{1}\Big(\sum_{j=1}^{m}p_{ijb}x_{ij}^{b} - 1 \Big) + \lambda_{2}\Big(\sum_{j=1}^{m}p_{ijb}x_{ij}^{b} - 1\Big)^{2},
    \label{QUBOconsistencyOneAsset}
\end{equation}
for all $i\in \mathcal{I}$. Essentially, the first of these terms favor solutions which satisfy the constraint whilst being as far away from the upper bound as possible.  The second term, instead, favors solutions which are as close to this upper bound. Effective tuning of the parameters is therefore necessary to balance the effects of each term. Note, however, in this case we do not need to relax the exposure constraint to an equality as we can better manage how far beyond the solution is from these lower bounds.
This equation is valid only for one asset, therefore we can modify equation (\ref{QUBOconsistencyOneAsset}) to consider all assets:
\begin{equation}
\begin{split}
   \lambda_{1}\sum_{i=1}^{n}\Big(\sum_{j=1}^{m}p_{ijb}x_{ij}^{b} - 1 \Big) \\
   + \lambda_{2}\sum_{i=1}^{n} \Big(\sum_{j=1}^{m}p_{ijb}x_{ij}^{b} - 1\Big)^{2}.
    \label{QUBOconsistencyConstraint}
\end{split}
\end{equation}

Following the same methods and the discussion on the one-to-one constraints used in the previous formulation, we can promote the remaining constraints to penalty terms in the QUBO.




The final QUBO can be written as:
\begin{equation}
\begin{split}
   & \lambda_{0}\sum_{i=1}^{n}\sum_{j=1}^{m} 
    {\Omega}_{ij}p_{ijb}x_{ij}^{b}\\
   & + \lambda_{1}\sum_{i=1}^{n}\Big(\sum_{j=1}^{m}p_{ijb}x_{ij}^{b} - 1 \Big) \\
   & + \lambda_{2}\sum_{i=1}^{n} \Big(\sum_{j=1}^{m}p_{ijb}x_{ij}^{b} - 1\Big)^{2} \\
    &- \lambda_{3}\sum_{j = 1}^{m}\Big(\sum_{i = 1}^{n}p_{ijb}x_{ij}^{b} a_{i}v_{i}H_{ij} - c_{j}\Big) \\
   &+  \lambda_{4}\sum_{j = 1}^{m}\Big(\sum_{i = 1}^{n}p_{ijb}x_{ij}^{b} a_{i}v_{i}H_{ij} - c_{j}\Big)^{2} \\
   & + \lambda_{5}\sum_{j=1}^{m}\sum_{g=1}^{G}\Big(\sum_{i=1}^{n}p_{ijb}x_{ij}^{b} T_{ig}a_{i} - K_{gj}\Big) \\
   & + \lambda_{6}\sum_{j=1}^{m}\sum_{g=1}^{G}\Big(\sum_{i=1}^{n}p_{ijb}x_{ij}^{b} T_{ig}a_{i} - K_{gj}\Big)^{2} \\
\end{split}
\end{equation}

\section{Numerical Illustrations}\label{sec:illustrations}

 \begin{table*}[ht!]
    \centering
    \begin{tabular}{cccccccc}\toprule[2pt]
         Formulation & Solver & \multicolumn{1}{c}{Cost Function} & \multicolumn{2}{c}{Consistency} & \multicolumn{2}{c}{Exposure}  & Objective Value\\ 
         &  &  $\lambda_0$& $\lambda_1$& $\lambda_2$& $\lambda_3$& $\lambda_4$ & \\
         \Xhline{0.02\arrayrulewidth}
        Balanced & D-Wave's Simulated Annealing Sampler & $10^3$ & $1$ & - & 1 & - & 0.5898 \\ \Xhline{0.02\arrayrulewidth}
        Balanced & Fujitsu's Digital Annealer & $10^5$ & $1$ & - & $300 $ & - & 0.7559  \\   \Xhline{0.02\arrayrulewidth}
        Unbalanced & D-Wave's Simulated Annealing Sampler & $1.5 \times 10^{4}$ & $1$ & $1$ & $1$ & $50$ & 0.5244   \\ \Xhline{0.02\arrayrulewidth}
        Unbalanced & Fujitsu's Digital Annealer & $2 \times 10^4$ & $1$ & $1$ & 1 & $50$ & 0.5803  \\ \Xhline{0.02\arrayrulewidth}
        Continuous LP & HiGHS (Simplex) & - & - & - & - & - & 0.4746 \\
        \Xhline{0.02\arrayrulewidth}
         
    \end{tabular}
    \caption{Values used for tuning the Lagrangian multipliers for each term in the QUBO for differing implementations. The objective value obtained for running this solution is also documented and for comparison, the value obtained by a continuous solver is displayed.}
    \label{tab:penalties}
\end{table*}

Utilizing the formulations presented in Sec. \ref{sec:level2} we now define a small instance of the \CO~problem based on a synthetic (however realistic) small dataset. We perform our tests on an Apple MacBook Pro with M2 Max processor and 16GB of memory. We remark that unlike the \KP~which even with a small instance of a few items has a relatively simple structure with a single constraint, the small \CO~problem instance we describe below has more complex constraints, and the interdependencies between the accounts and assets lead to a larger QUBO matrix and a much more challenging optimization landscape. 
In what follows, we utilize simulated annealing (SA), which as a metaheuristic algorithm, \cite{8130686} is quite sensitive to the problem structure and its performance can vary significantly depending on the problem instance. The increased complexity in the collateral optimization problem may make it harder for the SA to explore the solution space effectively, leading to suboptimal results or longer convergence times. In such cases, it might be beneficial to fine-tune the parameters of SA to improve its performance on more complicated problem instances. An interesting approach towards that direction was followed in \cite{sakuler2023real} which, for our context, we leave for future work.

Concretely, we have {\color{black} a portfolio of ten assets that have an approximate combined value of $ \$8.86$M. These assets can be categorized by their tier rating, $\omega = \{0.2, 0.5,0.8 \}$, into low, mid, and high tiered assets, respectively. Furthermore, the number of assets belonging to each category is chosen to be $4,2 $ and $4$, respectively. These assets are to be distributed in order to meet the requirements of 5 accounts. These requirements are distinguished by their duration, two of which are long-term and have a combined exposure of $ \sim \$1.49$M. The remaining are short-term requirements with a total exposure of $\sim \$1.09$M.

Due to restrictions, we slightly relax the problem by removing many-to-one constraints \eqref{eq:con3}. It is clear that in the absence of these constraints, this instance of the \CO~has a global optimum that can be easily obtained using classical strategies. As mentioned earlier, a compromise was needed between the precision of our results and the run-time performance, along with the limitation of the total number of bits that can be implemented in the solvers. To do this, we set the length of the bitstring representing each allocation, \eqref{eq:bitstring}, to be $7$. Hence, the granularity of our allocations is $1/127 \approx 0.0079 \%$,  (Since $M = 127$). In other words, the lowest percentage of available asset quantity we may post to any one account is approximately $0.0079\%$. This becomes important if an asset quantity is significantly large, as the corresponding solution will allocate more than what is necessary to meet the requirements. Also, if very strict limits were considered, many violations could occur if the bitstring length chosen was inadequate. To this effect, we ensure that our sample contains sensible values for the quantity of each asset so that they can be distributed with enough precision to satisfy the exposure requirements efficiently. } If no inequality constraints are included in the instance, then a total of $350$ qubits are required for both formulations ($10$ assets, $5$ accounts, $7$). However, due to the way we introduce the constraints, as discussed in the preceding section, we can reduce this number to $228$ qubits in the case of the unbalanced formulation and $298$ qubits for the balanced form, as the consistency constraint in this case still requires slack variables.

 Solving QUBO equations to obtain results that accurately reflect the goal of the objective function whilst simultaneously satisfying all constraints relies on the fine-tuning of the Lagrange multipliers. A potential solution is to make use of ``hybrid solvers'' such as D-Wave's Constrained Quadratic Model (CQM) \cite{d-wave} which automatically calculates these, we instead use an intuitive approach. The consistency constraint \eqref{eq:con1} is a hard constraint, as a solution that violates this constraint does not translate into a sensible business solution. Conversely, we treat the exposure requirement \eqref{eq:con4} as a soft constraint and allow for small violations to some margin $\epsilon$.  Additionally, for both balanced and unbalanced formulations, it is important to note that there are significant differences between the magnitudes of the coefficients of each term in the QUBO. This can make fine-tuning of the penalty weights a difficult task.  To manage this, we normalize each term in the QUBOs by dividing by their largest coefficient and then scale such that the lowest coefficient in each term has an order of magnitude of $1$. We chose the weight for the cost function to be a magnitude larger than those of the constraints, so that we achieve high-quality solutions. We then retrospectively increase the weight of the exposure and consistency terms to ensure that the system's constraints are satisfied.

 


\begin{figure}[!htb]
    \centering
    \includegraphics[scale=0.4]{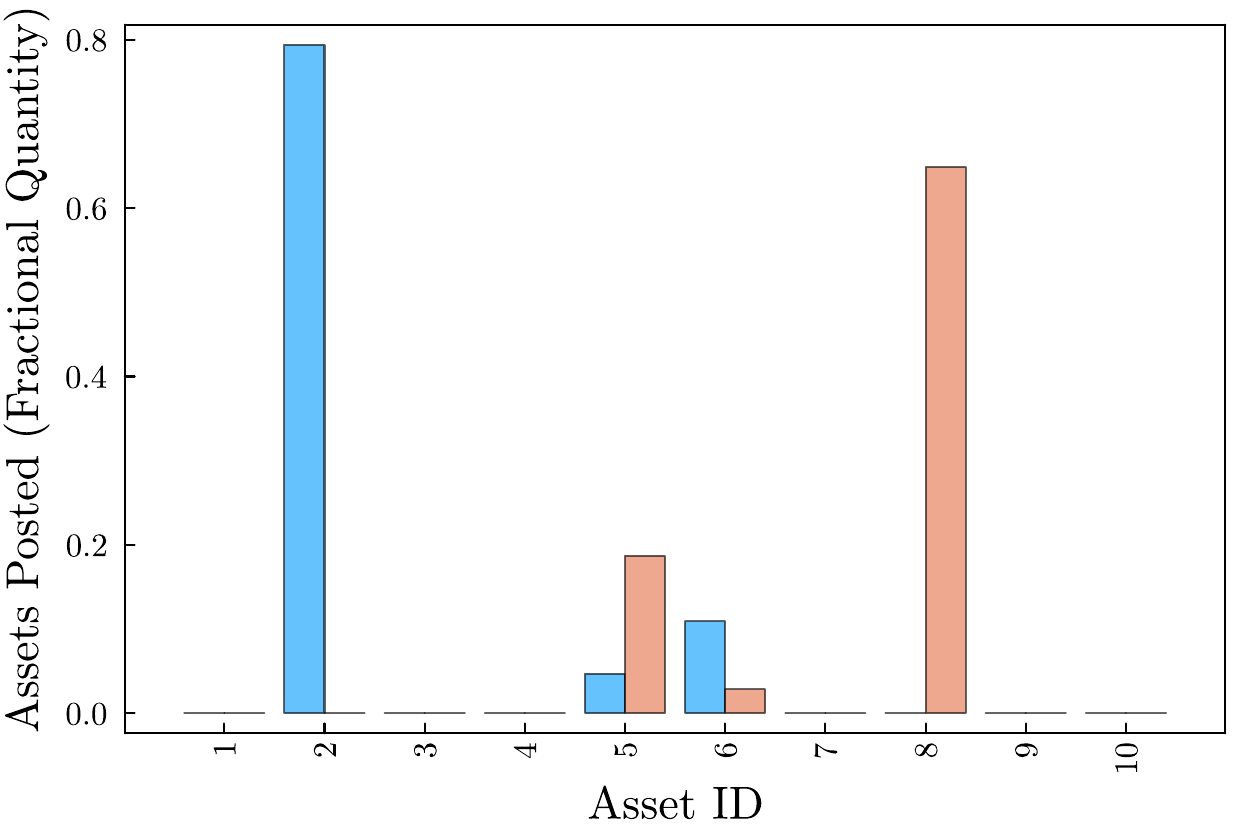}
    \caption{The optimal allocations of assets among accounts with short term (red) and long term (blue) requirements, determined through solving the \CO~ instance as a continuous LP with HiGHS. Assets IDs of $1-4$ are low-tier, $5-6$ are mid-tier and the final $7-10$ are the high-tier assets. }
    \label{fig:results_optimal}
\end{figure}

\begin{figure*}[!htb]
\centering
\subfigure{\includegraphics[scale=0.32]{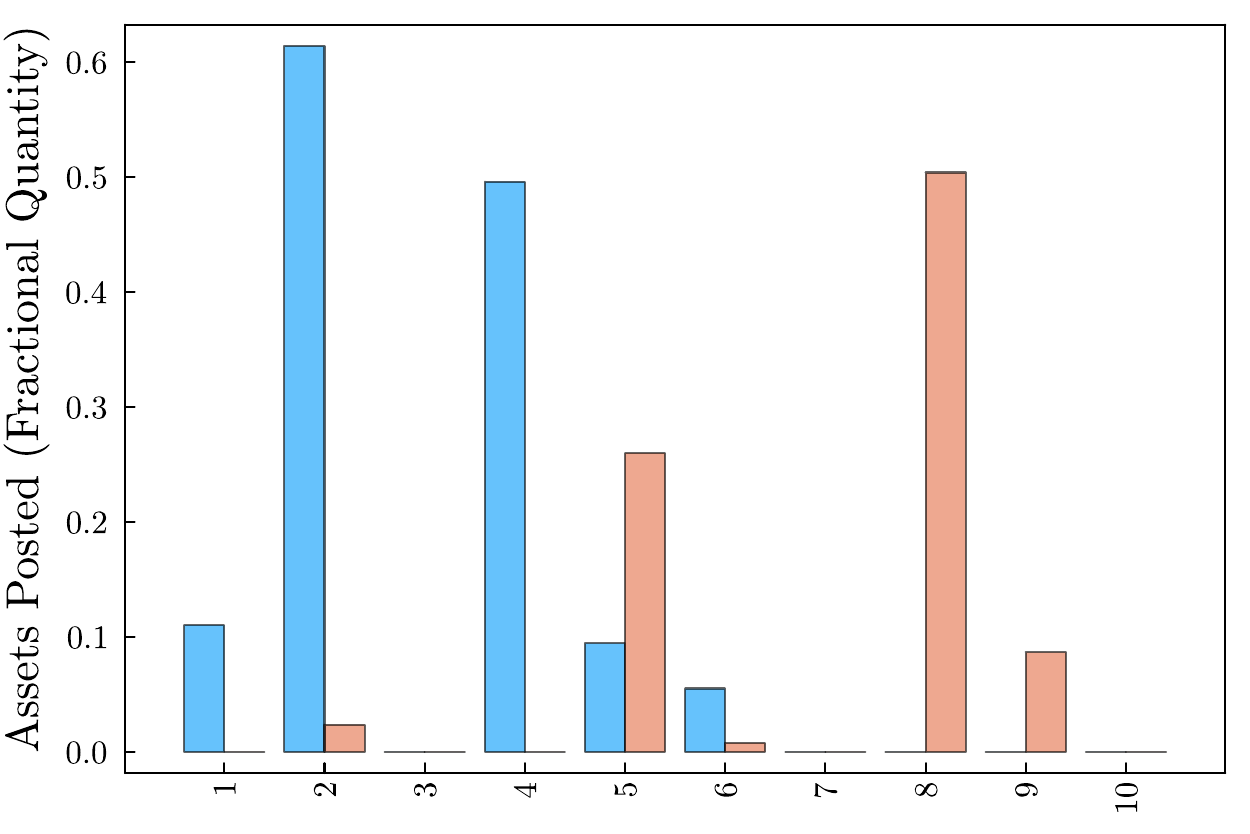}} 
\subfigure{\includegraphics[scale=0.32]{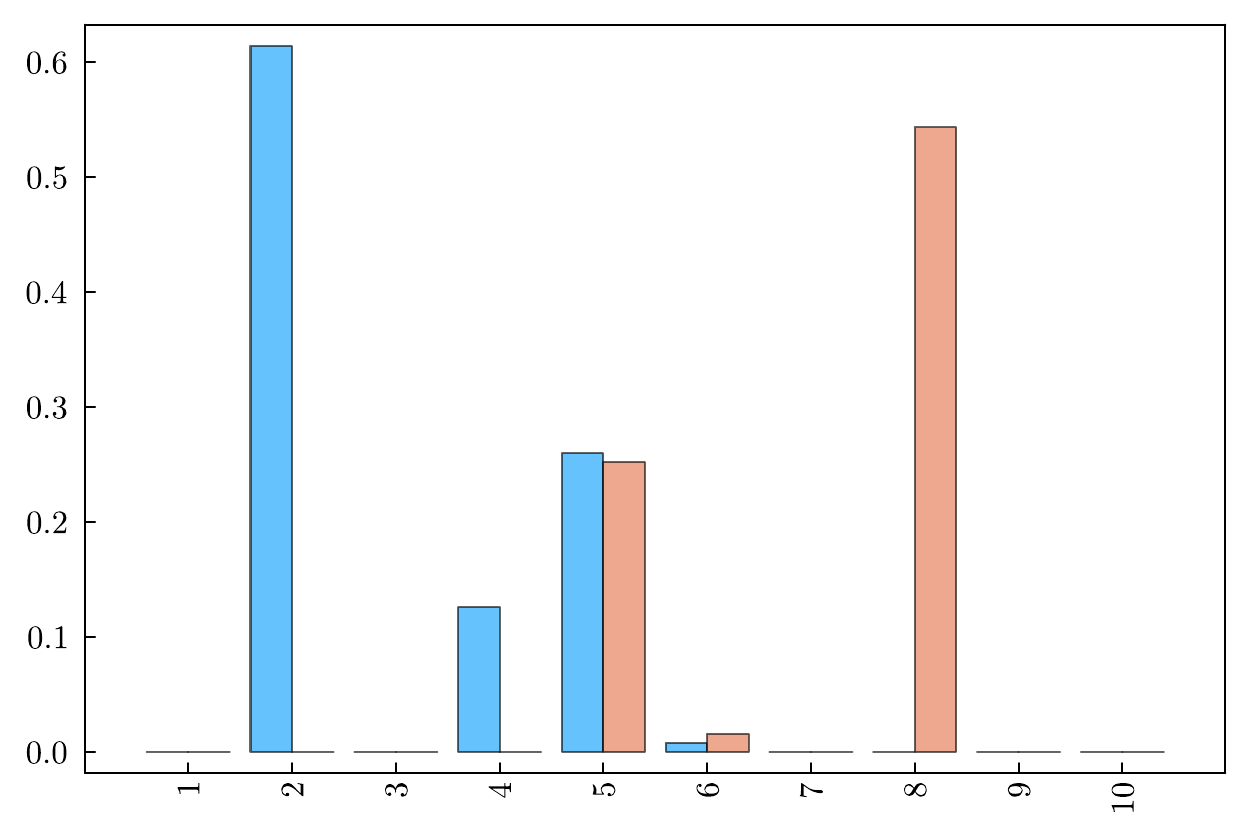}}
\subfigure{\includegraphics[scale=0.32]{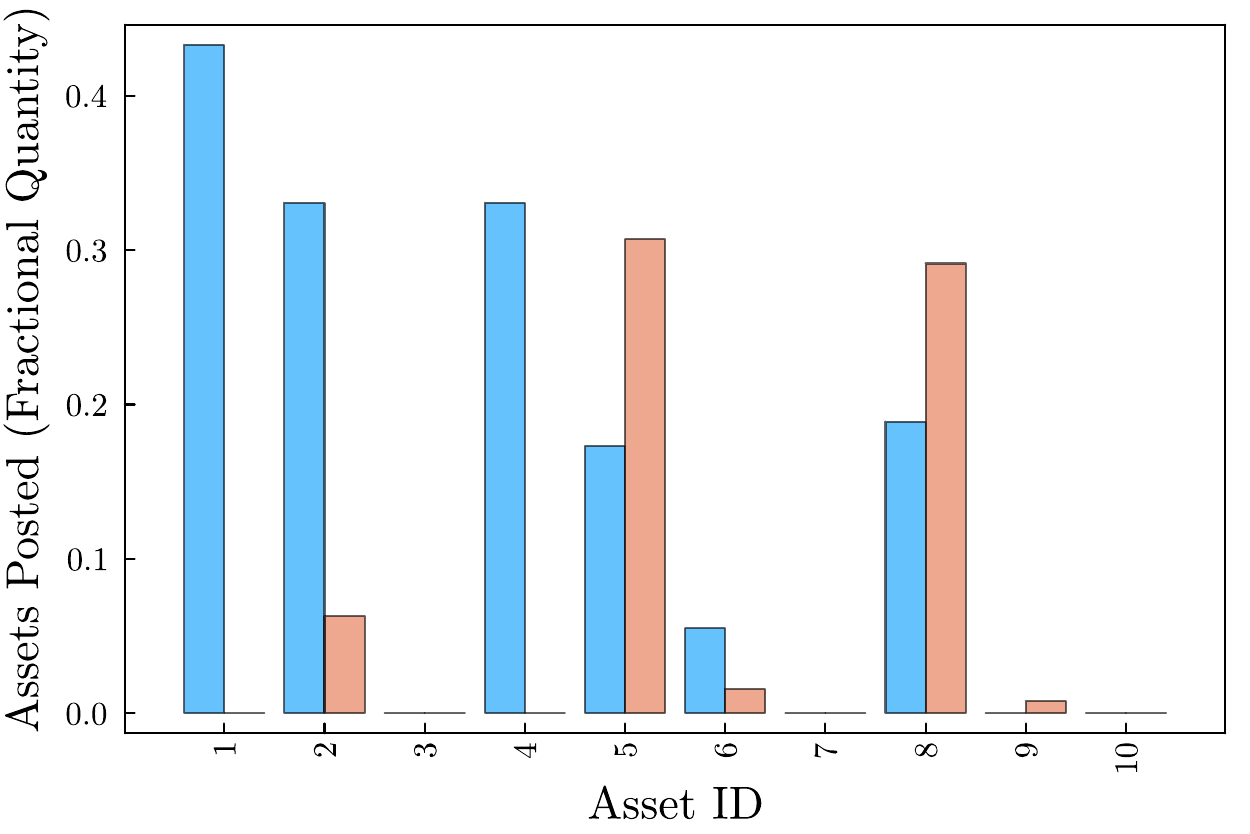}}
\subfigure{\includegraphics[scale=0.32]{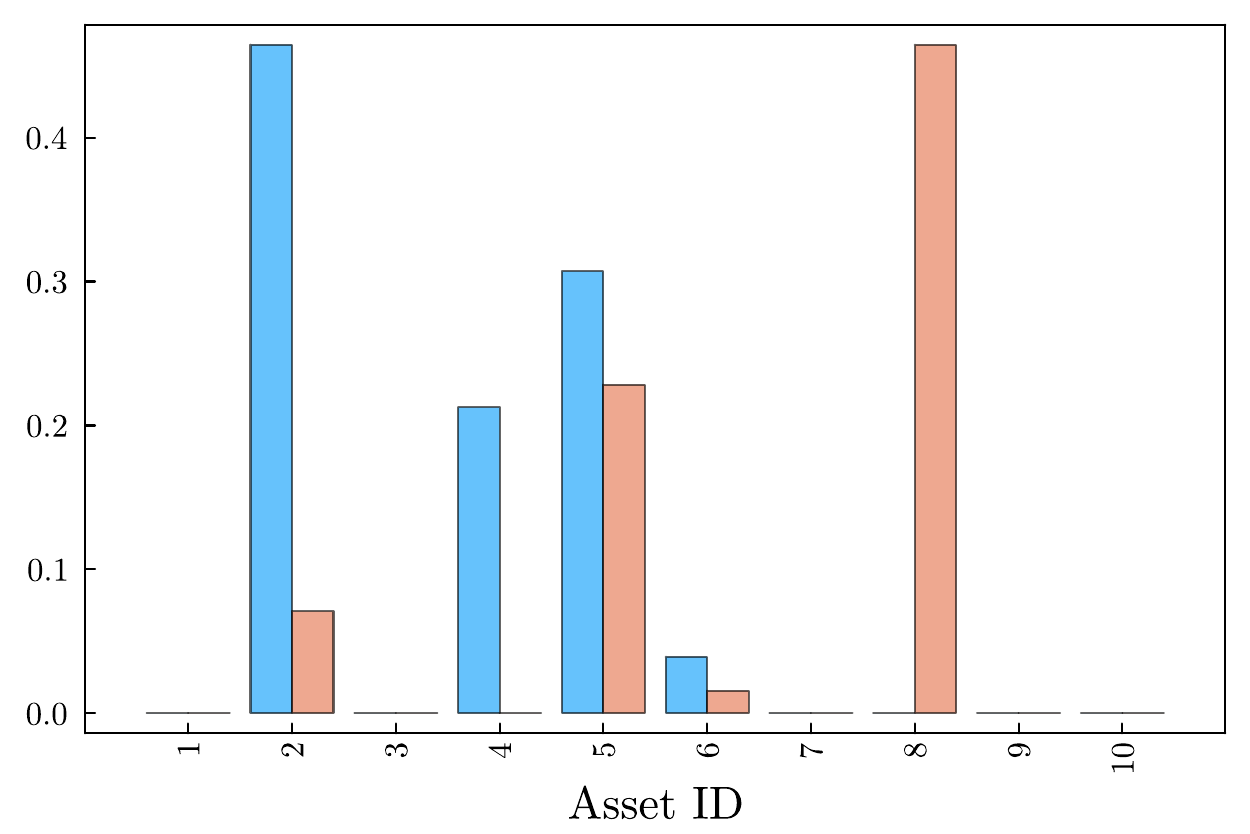}}
 \caption{The allocation, of different assets amongst accounts with short term (red) and long term (blue) requirements. Results are determined by D-Wave's simulated annealer (top) and Fujitsu's Digital annealer (bottom) with balanced (left) and unbalanced (right) formulations. The asset IDs of $1-4$ are low-tier, $5-6$ are mid-tier and the final $7-10$ are the high-tier assets.}
 \label{fig:resultsa}
\end{figure*}

Overall, our results were mixed in the sense that none of our runs managed to reach the global optimum nor to produce the globally optimal allocation with each run converging in a different local minimum. We estimate that one reason for this behavior is the limited number of runs performed that do not allow the annealing process to explore sufficient search space. This can be easily redeemed by increasing the number of runs (potentially decreasing the step size) and utilizing more compute power. However due to limited resources and compromising between computational runtime and the accuracy of the solution, we opted to chose a modest number of runs.

\begin{figure}
    \centering
    \includegraphics[scale=0.40]{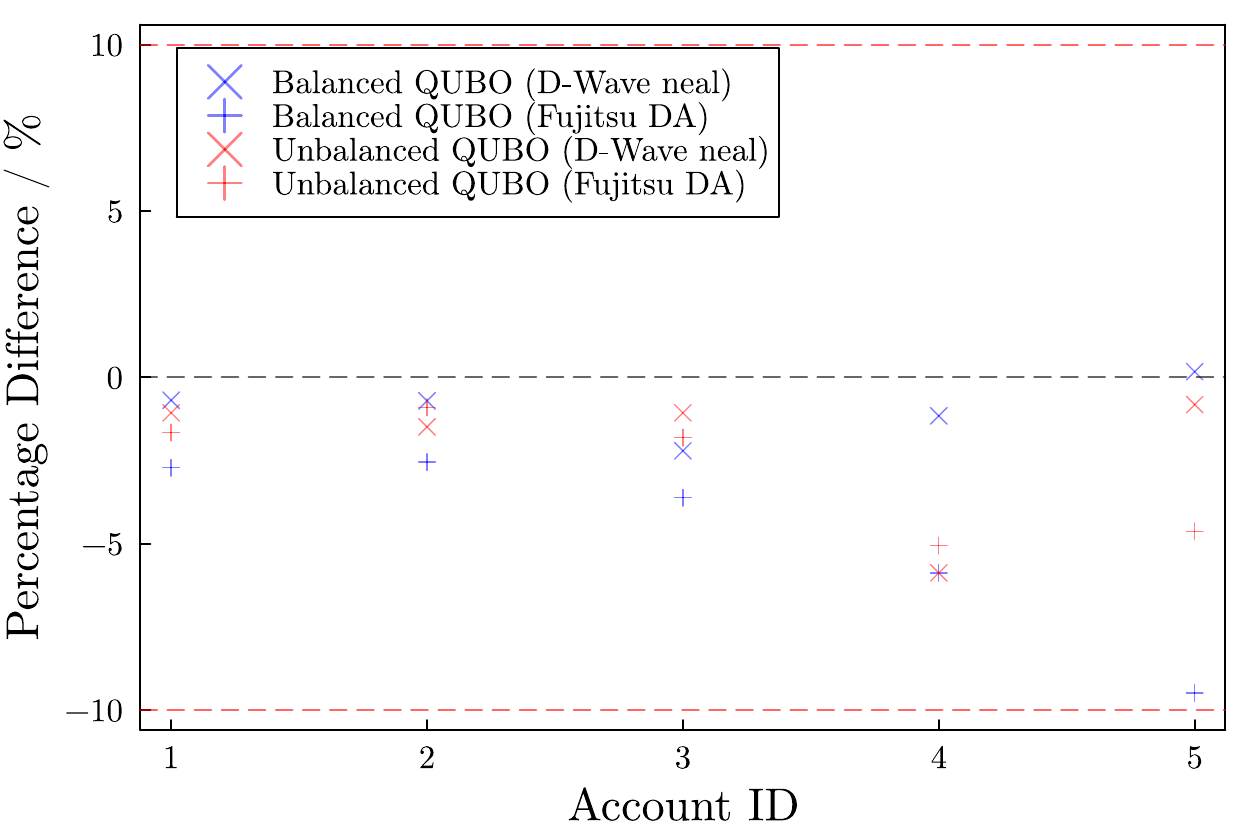}
    \caption{The percentage of the exposure requirements that have been met for each account. The dashed line represents the solution given by HiGHS solver which perfectly meets each requirement. We see greater deviations for the requirement of account 4, which is due to its outstanding exposure being an order of magnitude less than those of the other accounts, hence it has a lower weighting in the QUBO.}
    \label{fig:Overhead_exposure}
\end{figure}

Table \ref{tab:penalties} displays the values chosen for each of the penalty weights and the resultant objective value that was outputted. Fig. \ref{fig:results_optimal} shows the global optimal solution solved using HiGHS. The asset allocations using  {neal} and Fujitu's Digital Annealer are shown in Fig. \ref{fig:resultsa}. Additionally, Fig. \ref{fig:Overhead_exposure} displays the percentage differences, for each of the solver, between the total values posted and the required exposures for all the accounts.

 



\section{Summary and conclusions}\label{sec:summary}

We surveyed the problem of collateral optimization from a business perspective and provided a business-realistic MILP formulation suitable to be mapped to a QUBO. In turn, we provided two QUBO formulations, one based on slack variables and one based on unbalanced penalization, and we implemented a small \CO~problem instance on small emulations of hybrid solvers. We observed that the unbalanced penalization approach yields objective function values much closer to the global optimum, obtained using the simplex method, while the slack-based (balanced) approaches were further off. To this end, both approaches fail to find the global optimum even for such a relatively small, but non-trivial, problem instance. While we did not run our computations on specialized hardware (quantum annealer or digital annealer), to some extent, we indeed rediscover the expectation that heuristic approaches fail to find globally optimal solutions very often, including when implemented in quantum hardware \cite{vert2021benchmarking}. The ``take-home message'' from the numerical illustrations clearly showcases the relevance of the formulation chosen for a given problem when one is interested in making it quantum-ready.

Note, however, that the aim of this paper is to show the relevance as well as the technical formulation of the collateral optimization problem for quantum or hybrid solutions rather than to perform a detailed benchmarking and, as a matter of fact, several improvements can be performed in order to obtain higher quality solutions (warm starting techniques, optimizing the annealing schedule, QUBO parameter optimization, utilization of GPU and tensor cores, and explore improved methods for SA \cite{nakano2016improved}). 

The question of whether MILPs such as the collateral optimization problem should be cast as QUBOs and approached by heuristic solvers (quantum or not) is tricky. To some extent, classical solvers such as Gurobi, CPLEX, or even HiGHS perform exceptionally well in many large instances. Improving thereof could be less or more beneficial depending on the situation or problem instance and its complexity. If a ``quantum'' solution is preferred, finding the ultimate strategy for this approach is a crucial problem (formulation of the problem, hyperparameter choice, solver choice, etc).

The \CO~problem, as presented above, can be extended in a variety of ways. The formulations provided above are in no way unique and further investigation might yield better results even for the local solvers. For example, recent research \cite{PRXQuantum.3.030304} has shown that an alternative formulation based on implicit penalization by restricting the Hamiltonian dynamics, in the context of parallelized QAOA, can be powerful and such an approach can be suitable for the collateral optimization problem as well. A different approach based on a stochastic quantum Monte Carlo algorithm, that mimics quantum annealing, was proposed \cite{https://doi.org/10.48550/arxiv.2302.12454} for such large-scale optimization problems making it another suitable candidate for \CO. The benefit of this approach is that it can handle fully connected graphs which can be the case in certain \CO~instances.
\vspace{1em}




\noindent
\emph{Acknowledgements}.
We would like to thank Jakub Marecek, Vyacheslav Kungurtsev, David Snelling, Reiss Pikett and Philip Intallura for useful discussions and suggestions. 
\vspace{1em}

\noindent
\emph{Disclaimer}.
This paper was prepared for information purposes
and is not a product of HSBC Bank Plc. or its affiliates.
Neither HSBC Bank Plc. nor any of its affiliates make
any explicit or implied representation or warranty and
none of them accept any liability in connection with
this paper, including, but limited to, the completeness,
accuracy, reliability of information contained herein and
the potential legal, compliance, tax or accounting effects
thereof. Copyright HSBC Group 2023.

\noindent

\bibliographystyle{unsrt}
\bibliography{main.bib}


\begin{appendices}

\section{Collateral optimization terminology }\label{ap:explanation}
Below, in Table \ref{tab:TableofDefinitions}, we collect a few of the common terms used in the context of \CO~and, more generally, collateral management.

\begin{table}[hbt!]
\begin{center}
\begin{tabular}{p{0.1\textwidth} p{0.37\textwidth}}
\toprule[2pt]
 {\bf Term} & {\bf Definition} \\ 
 \hline
  {\bf Collateral} & Any valuable asset used within lending agreements which can be seized by the lender from the borrower if they fail to repay the loan. This can consist of any asset deemed valuable by the lender. Acceptable types of assets are: equity (stocks), bonds and cash. \\
  {\bf Market Value} & The price that an asset could be sold for when auctioned in the marketplace.\\
  {\bf Haircut} & A reduction applied to the market value of an asset. Used to ensure that the posted asset will cover the cost of the outstanding exposure, given the risk that an assets value can fluctuate. \\
  {\bf Exposure} & The amount of capital that one stands to lose (risk) if an investment fails. \\
 \hline
\end{tabular}
\end{center}
    \caption{Definition of financial terms used.}
    \label{tab:TableofDefinitions}
\end{table}

Table \ref{tab:TableOfNotation} collects all relevant \CO~related quantities/variables used in the main body of this paper.

\begin{table}[hbt!]
    \begin{center}
\begin{tabular}{p{0.1\textwidth}  p{0.35\textwidth}}
\toprule[2pt]
 {\bf Symbol} & {\bf Explanation} \\ 
 \hline
  $\mathcal{I}$ & Set of $n$ inventory assets. \\
  $i$ & Index representing a particular asset.  \\
  $a_i$ &  Maximum available quantity of asset $i$.\\
  $v_i$ & Dollar value of a single unit of asset $i$. \\
  $\omega_i$ & Tier value representing the quality of asset $i$. \\
  $\mathcal{A}$ & Set of $m$ accounts.\\
  $j$ & Index representing a particular account. \\
  $c_{j}$ & Dollar value exposure of account $j$ required. \\
  $d_{j}$ & Binary value indicating whether account obligation $j$ is short or long term. \\
  $\Omega_{ij}$ &  Tier value of asset $i$ depending duration of \\ & account $j$. \\ 
  $Q_{ij}$ & Decision variable indicating the percentage \\ & of asset $i$ allocated to account $j$.\\
  $B_{ij}$ & Limit on individual allocation of asset $i$ to \\ & account $j$.\\
  $H_{ij}$ & Haircut factor used to reduce the value of an asset. \\
  $K_{gj}$ & Limit on the allocation of assets in group $g$ to \\ & account $j$.\\
  $T_{ig}$ & Binary value indicating whether asset $i$ \\ & belongs to group $g$.\\
 \hline
\end{tabular}
\end{center}
    \caption{List of all symbols used in the formulation of the collateral optimization problem. \ref{ap:explanation}.}
    \label{tab:TableOfNotation}
\end{table}

\section{Constraint Penalization in QUBOs and the Ising model}\label{ap:qubo_to_ising}

In the main body of this paper, we have seen that in order to encode penalty terms into QUBOs we need to use, for example, slack variables. Here we summarize how we proceed to do so. In the \KP~as well as the collateral optimization problem, we have constraints of the form 
\begin{align}
    Ax &\leq c \\
    Bx &\geq d. \
\end{align}
Let us consider the former constraint since the latter is essentially the same up to a factor of $-1$. It is common to view this constraint as $Ax-c \leq 0$ and then introduce a slack variable $s\in \mathbb{R}_{\geq 0}$ such that we define the penalty term $Ax-c + s = 0$. This is defined such that 
\begin{align}
    \sup s =  c - \inf Ax.
 \end{align}

Then, as mentioned already in Sec. \ref{sec:knapsack}, one needs to map a QUBO to the Ising Hamiltonian that performs Quantum Annealing (QA), a restricted model of adiabatic quantum computation. Let us now briefly explain this connection. 

The (classical) Ising model was first introduced as a mathematical model of ferromagnetism \cite{cipra1987introduction}. The variables take values in a discrete (binary) set $ \sigma_i=\{\pm 1\}$, and are typically referred to as \emph{spin}, since in the physical model they describe the atomic spin of the particles. The model consists of a lattice, $\Lambda$, with each lattice site, $j\in\Lambda$, having an assigned spin $\sigma_j$. The energy of a specific spin configuration is measured by the Hamiltonian of the system
\begin{equation*}
    H= -\sum_{j<k}J_{j,j+1}\sigma_j\sigma_{j+1}-\sum_{j} h_j\sigma_j + \varepsilon,
\end{equation*}
where $J_{j,j+1}$ encodes the ``nearest neighbour interaction'' between adjacent sites, and $h_j$ describes some external field that interacts with each site \cite{lucas2014ising}. The overall minus sign is just a convention\footnote{However, the sign of $J_{jk}$ does have a more interesting interpretation. If $J_{jk}>0$ the model describes ferromagnetism, while if $J_{jk}<0$ it describes anti-ferromagnetism. } and $\varepsilon$ refers to a constant energy overhead.

In the quantum version of the Ising Hamiltonian, the sites are represented by a qubit. The spins are then simply given by the Pauli matrix $\sigma_{j}^z$ acting on the $j$'th site, with eigenvalues $\pm1$ when acting on the computational basis states $\{|0\rangle, |1\rangle\}$. In physics, the interesting problem is typically to find the ground state energy, or lowest energy eigenstate, of the Hamiltonian.

One thing we can immediately notice is that this model is quadratic in the spins and thus resembles the QUBO. By the simple change of variables $x_j=\tfrac{1-\sigma_j}{2}$ we directly see that the spin variables $\pm1$ are mapped to the binary variables $\{0,1\}$. In the quantum version we map the QUBO variables $x_j$ to the operators $x_j=\pmb{1}\otimes...\otimes\tfrac12(\pmb{1}-\sigma_z)\otimes\pmb{1}\otimes...,$ where the non-trivial operator acts on the $j$'th site. This operator has eigenvalues $0$ or $1$ when acting on the states $|0\rangle$ respectively $|1\rangle$. Trough this change of variables, the QUBO problem is thus equivalent to finding the ground state energy of the Ising Hamiltonian. When mapping between the QUBO and the Ising model we might also need to account for the fact that we only minimise over the Ising model, while sometimes, for example in the \KP, we are seeking a maximum. This is of course easily accounted for by changing the relevant signs.

As an example, the \KP~\eqref{eq:KPslack1} has an Ising Hamiltonian with:
\begin{align*}
    J_{j,j+1} = -\frac{\lambda_{0}}{4}w_{j}w_{j+1}, \\
    h_{j} = -\left(\frac{v_{j}}{2} + \lambda_{0}Kw_{j}\right), \\
    \varepsilon = -\frac{\lambda_{0}}{4}\sum^{n}_{i}w_{i}^{2} - \frac{1}{2}\sum^{n}_{i}v_{i} -\lambda_{0}K^{2},
\end{align*}
where we define $K = S - W + \frac{1}{2}\sum^{n}_{i=1}w_{i}$ for convenience. 

\end{appendices}

\EOD

\end{document}